\documentclass[a4paper,DIV=11,abstracton]{scrartcl}

\usepackage[T1]{fontenc}
\usepackage[english]{babel}
\usepackage[bf,normal]{caption}
\usepackage{amsmath,amsfonts,amsthm,mathtools}
\usepackage{url}
\usepackage{graphicx}
\usepackage{todonotes}
\usepackage{appendix}
\usepackage{caption,subcaption}
\usepackage{enumitem}
\usepackage{booktabs}
\usepackage[ruled,vlined,linesnumbered]{algorithm2e}
\usepackage{xcolor}
\usepackage{authblk}
\usepackage[pdffitwindow=false,
            plainpages=false,
            pdfpagelabels=true,
            pdfpagemode=UseOutlines,
            pdfpagelayout=OneColumn, 
            bookmarks=true,
            bookmarksopen=true,
            colorlinks=true,
            hyperfootnotes=false,
            linkcolor=blue!30!black,
            urlcolor=blue!30!black,
            citecolor=green!50!black]{hyperref}

\newcommand{\R}{\mathbb{R}}
\newcommand*{\EE}{\mathbb{E}}

\theoremstyle{plain} % italic
\newtheorem{thm}{Theorem}[section]
\theoremstyle{definition} % not italic
\newtheorem{remark}[thm]{Remark}

\author[1,2]{Jan-Hendrik Niemann}
\author[1]{Samuel Uram}
\author[1]{Sarah Wolf}
\author[2]{Nata\v{s}a Djurdjevac Conrad}
\author[2]{Martin Weiser}
\affil[1]{Department of Mathematics and Computer Science, Freie Universit\"at Berlin, Germany,}
\affil[2]{Zuse Institute Berlin, Germany\footnote{This research has been funded by the Deutsche Forschungsgemeinschaft (DFG, German Research Foundation) under Germany's Excellence Strategy MATH\texttt{+}: The Berlin Mathematics Research Center (EXC-2046/1, project ID: 390685689). Email: \texttt{\{niemann,natasa.conrad,weiser\}@zib.de}, \texttt{sarah.wolf@fu-berlin.de}, \texttt{usamuelsunadh@gmail.com}}}

\title{Multilevel Optimization for Policy Design with Agent-Based Epidemic Models}
\date{March 22, 2023}

\begin{document}

\maketitle

\begin{abstract}
    Epidemiological models can not only be used to forecast the course of a pandemic like COVID-19, but also to propose and design non-phar\-ma\-ceu\-ti\-cal interventions such as school and work closing. In general, the design of optimal policies leads to nonlinear optimization problems that can be solved by numerical algorithms. Epidemiological models come in different complexities, ranging from systems of simple ordinary differential equations (ODEs) to complex agent-based models (ABMs). The former allow a fast and straightforward optimization, but are limited in accuracy, detail, and parameterization, while the latter can resolve spreading processes in detail, but are extremely expensive to optimize. We consider policy optimization in a prototypical situation modeled as both ODE and ABM, review numerical optimization approaches, and propose a heterogeneous multilevel approach based on combining a fine-resolution ABM and  a coarse ODE model. Numerical experiments, in particular with respect to convergence speed, are given for illustrative examples.
\end{abstract}


\section{Introduction}

The global COVID-19 pandemic has highlighted as never before the need for mathematical modeling to forecast infection spreading and assess consequences of various non-phar\-ma\-ceu\-ti\-cal interventions taken to counter it, such as closing of schools, social distancing rules, or travel restrictions. Such intervention policies have different direct costs as well as social or economic impacts, and differ in their effectiveness with respect to various objectives (e.g.,~minimizing the expected number of deaths or the number of new infections in a given time window), and in their negative side effects. While policy making must take into account many quantitatively poorly understood aspects, and therefore cannot rely solely on mathematical models of infection spreading, model-based policy optimization can be an important decision support tool. For example, optimization may reveal policies or combinations of policies that were not obvious from common sense alone, or that may avoid assumed trade-offs or produce unexpected benefits.

Models used to study epidemics range from deterministic compartmental models based on deterministic ordinary differential equations (ODEs) such as the well-known Kermack--McKendrick SIR model (considering population fractions of susceptible, infected, and recovered individuals) to spatially resolved stochastic agent-based models (ABMs)~\cite{MaratheRamakrishnan2013, Brauer2017, TracyCerdaKeyes2018}, which describe the dynamics of infection as arising from a large number of discrete interactions between agents in complex interaction networks embedded in a common environment. 

Each modeling approach has its own trade-offs. On the one hand, ODE models are computationally inexpensive and easy to understand and analyze, so a variety of mathematical methods and tools are readily available for their analysis. On the other hand, as aggregate models, ODE models provide coarse rather than detailed forecasts of infection dynamics that do not allow for uncertainty estimates of the forecasts, and policy measures are difficult to translate into changes in model parameterization. Models based on stochastic differential equations (SDEs) can serve as a means to quantify uncertainty. However, these are often still compartmental models with the same drawbacks.

In contrast, complex ABMs allow for the representation of locally adapted policies and the analysis of their consequences in terms of spatial (e.g., neighborhood, city, region), compartmental (e.g., children, students, workers, retirees), or sectoral (e.g., medical, educational, retail) dimensions, and thus provide much more detailed forecasts, including uncertainty quantification. This makes ABMs attractive for policy design, especially for locally adapted policy design in complex environments. The main drawback is that it is computationally expensive to simulate realistic ABMs. Also, the models are inherently stochastic, requiring the simulation of several to many samples to obtain a reliable answer. Perhaps most importantly, mathematical formulation (and hence analysis) is still the exception (for some examples see~\cite{izquierdo+al2009, hinkelmann+al2011, banisch+al2012, Conrad2018EPJ, HelfmannDjConrad2019, niemann+al2021}) rather than the rule. While ABMs are widely used in various scientific disciplines such as sociology, economics, or geography, the focus is mostly on computer simulations.

Optimal control is about finding a control (usually a time-dependent function) for a dynamic system and a given time window in such a way that an objective function is optimized, often under given constraints. Numerical methods for solving complex and nonlinear optimization problems are well established. Optimal control and stability analysis of deterministic ODE epidemic models to support policy design has a long tradition~\cite{Wickwire1977}, but is still an active area of research, e.g., for tuberculosis~\cite{MoualeuWeiserEhrigDeuflhard2015} or HIV~\cite{DuwalWinkelmannSchuetteKleist2015}. Most recently, optimal control has been used to study COVID-19. In~\cite{Olivares2022}, robust optimal control and polynomial chaos expansion are used to find optimal vaccination strategies under uncertainty about compliance with social distancing and testing strategies. In further works, e.g.,~\cite{Ullah2020,Charpentier2020,Lazebnik2021}, the effects of non-phar\-ma\-ceu\-ti\-cal interventions such as different lockdown scenarios are studied using optimal control.

In contrast, mathematical optimization of ABMs faces significant challenges. Basic concepts necessary for optimization, such as derivatives, are not well-defined due to the discrete nature of the steps and decisions made by agents in most ABMs. A recent idea to circumvent this is internal smoothing~\cite{andelfinger2021}. Adjoint concepts for discontinuous option pricing~\cite{Fries2021} may also provide a direction for future research. However, most of the approaches used so far avoid the consideration of derivatives in the ABM altogether, for example, by applying genetic algorithms directly to an ABM~\cite{schutte2010,Miikkulainen2021}, or by relegating the optimization problem to another model that approximates the ABM. This may be a reduced ABM, derived for example by coarsening the spatial resolution of the underlying one~\cite{oremland+laubenbacher2014}, so that again heuristic algorithms are applied for optimization. In most cases, an ``equation-based'' or ``system-level'' model is used, for which optimization methods exist. The choice of such a coarser model is based on what is already available for similar phenomena, including difference or ordinary and partial differential equations. Alternatively, a simple functional form is chosen and parameters are fitted to the ABM output (see~\cite{an+al2017} and references therein). The optimal controls for the reduced model are then lifted back to the ABM. 

In this paper, we apply a multilevel approach that approximates the original ABM at the fine level with a deterministic ODE at the coarse level, using the well-known epidemic SIR model. Multilevel methods exploit the computationally cheaper optimization of coarser models to reduce the optimization effort in a more detailed ``fine'' model. However, the perspective taken is not that of finding optimal interventions for infectious diseases such as COVID-19, but that of numerical optimization, since finding optimal interventions is the subject of a political and societal debate that we do not intend to enter in this article. To minimize a given (prototypical and exemplary) objective, we consider derivatives and descent along them. In the (stochastic) ABM, where exact derivatives cannot be obtained due to its inherent discontinuity, we work with the numerical approximation of derivatives by finite differences. In the multilevel approach presented here, the combination of coarse and fine levels is iterative: a coarse-level ODE model is used to generate trial steps for descent directions in the fine-level ABM. This step must then pass an acceptance test to ensure descent of the objective with high confidence.

The remainder is organized as follows. In Section~\ref{sec:epidemic-models}, we introduce three epidemiological models and investigate their quantitative relation. In particular, we present a classical SIR-type ODE and an associated homogeneous ABM represented as a Markov jump process as stylized examples. As a third model, we briefly present a detailed, community-specific, and geographically referenced ABM developed by~\cite{goldenbogen+al2020} as a less stylized, more realistic example. Policy optimization is considered in Section~\ref{sec:optimization}, where we consider policy definition and realization, as well as the formulation of an optimization objective, in a simple, prototypical situation. We briefly review state-of-the-art optimization algorithms suitable for ODEs and ABMs and discuss their computational complexity. We then introduce a heterogeneous multilevel optimization approach that combines a fine-resolution ABM and a coarse ODE approximation to speed up convergence. These algorithms, and in particular their relative computational efficiency, are investigated in Section~\ref{sec:num-examples} using numerical examples of the transmission dynamics of SARS-CoV-2. Open questions and future work are discussed in Section~\ref{sec:conclusion}.

%==============================================================================================

\section{Epidemiological Models} 
\label{sec:epidemic-models}

This section introduces terms and notation for the type of epidemiological models we use, namely models based on deterministic ODEs in Section~\ref{sec:basicODE} and stochastic ABMs in Section~\ref{sec:ABM}. The relation and agreement of these types of models are discussed in Section~\ref{sec:agreement}.

\subsection{Deterministic ODE Models}
\label{sec:basicODE}

\begin{figure}[t]
    \centering
    \begin{subfigure}{0.49\textwidth}
        \centering
        \caption{}
        \includegraphics[width=\textwidth]{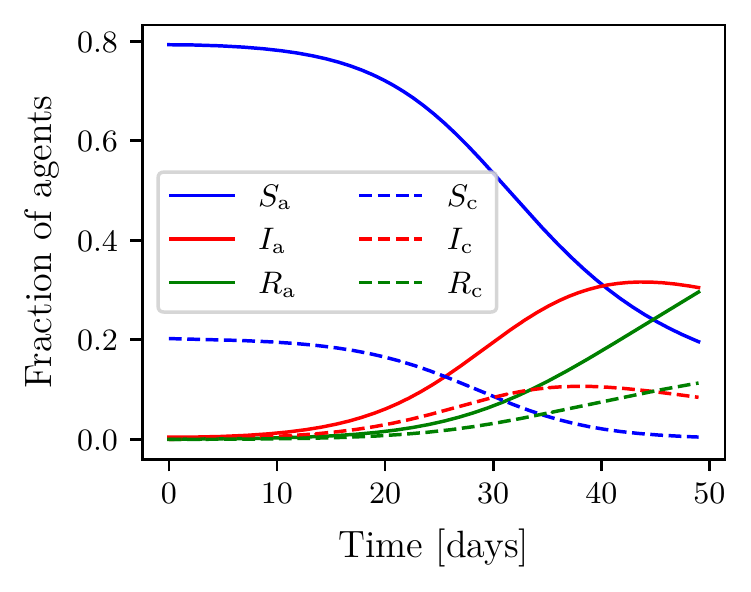}
    \end{subfigure}
    \begin{subfigure}{0.49\textwidth}
        \centering
        \caption{}
        \includegraphics[width=\textwidth]{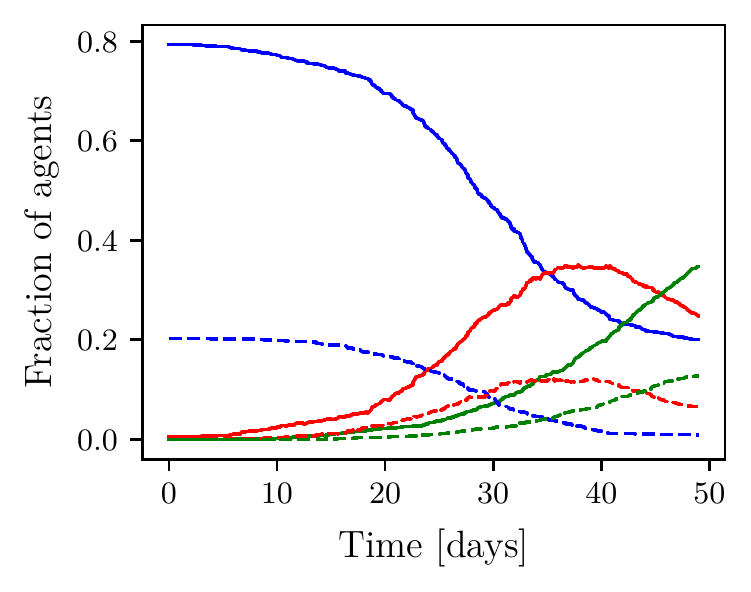}
    \end{subfigure}
    \caption{(a)~Numerical solution of the ODE~\eqref{eq:ODE_SIR} and (b)~aggregated trajectory of a single simulation of GERDA given identical initial conditions divided into adults (solid lines) and children (dashed lines) for parameters given in Table~\ref{tab:parameters}.}
    \label{fig:sample-ODE-and-ABM-run}
\end{figure}

Mathematical models based on ODEs are commonly used to simulate epidemic spread on a macro-scale of a fully mixed population, e.g., a city or a country. These models separate the population into different compartments (therefore often referred to as \emph{compartmental models}) and count only the proportions of individuals belonging to each compartment. Compartments can, for example, be based on the infection status or demographics or both. In this paper, we consider so-called SIR models~\cite{Brauer2017} and divide the population into two age groups: adults and children (compartments denoted by subscript $\mathrm{a}$ or $\mathrm{c}$, respectively). Separation into age groups has proven beneficial due to different infection rates, see, e.g.,~\cite{Cai2020} for the case of SARS-CoV-2. Individuals who have not yet been infected belong to the susceptible compartment $S_\circ$, for $\circ\in\{\mathrm{a},\mathrm{c}\}$. With a certain infection rate $r_{*\to \circ}$, for $*,\circ\in\{\mathrm{a},\mathrm{c}\}$, individuals from $S_\circ$ can get infected by an individual from $I_*$ and move to the compartment $I_\circ$. An infection is caused by interactions of the type $S_\circ + I_* \rightarrow I_\circ + I_* $, which are referred to as \emph{second-order} interactions since two individuals are involved. This notation stems from the chemical context, where it indicates a transition of two particles of the given types to two particles with the new types on the right-hand side of the arrow. Individuals in the group $I_\circ$  recover with rate $r_\circ$ and move to the compartment $R_\circ$ due to the \emph{first-order} reaction $I_\circ \rightarrow R_\circ$. The system of ODEs describing this model is given by
\begin{equation} \label{eq:ODE_SIR}
    \begin{aligned}
        \dot S_\mathrm{a} &= - S_\mathrm{a}\, (r_{\mathrm{a} \to \mathrm{a}} I_\mathrm{a} + r_{\mathrm{c} \to \mathrm{a}} I_\mathrm{c})\\
        \dot S_\mathrm{c} &= - S_\mathrm{c}\, (r_{\mathrm{c} \to \mathrm{c}} I_\mathrm{c} + r_{\mathrm{a} \to \mathrm{c}} I_\mathrm{a})\\
        \dot I_\mathrm{a} &= S_\mathrm{a}\, (r_{\mathrm{a} \to \mathrm{a}} I_\mathrm{a} + r_{\mathrm{c} \to \mathrm{a}} I_\mathrm{c}) - r_\mathrm{a} I_\mathrm{a}\\
        \dot I_\mathrm{c} &= S_\mathrm{c}\, (r_{\mathrm{c} \to \mathrm{c}} I_\mathrm{c} + r_{\mathrm{a} \to \mathrm{c}} I_\mathrm{a}) - r_\mathrm{c} I_\mathrm{c}\\
        \dot R_\mathrm{a} &= r_\mathrm{a} I_\mathrm{a} \\
        \dot R_\mathrm{c} &= r_\mathrm{c} I_\mathrm{c}
    \end{aligned}
\end{equation}
Obviously, the sum of the groups $S_\circ + I_\circ + R_\circ$ is conserved over time, such that $R_\circ$ is often omitted because of redundancy. If it is needed, it can be computed from $S_\circ$ and $I_\circ$ as $R_\circ = N_\circ - S_\circ - I_\circ$, where $N_\circ$ denotes the fraction of individuals of the respective age group. In the following, when we do not distinguish between age groups, we will refer to the number of susceptible $S = S_\mathrm{a}+S_\mathrm{c}$, infected $I = I_\mathrm{a}+I_\mathrm{c}$, and recovered individuals $R = R_\mathrm{a}+R_\mathrm{c}$, such that $S+I+R = 1$. A numerical solution of ODE~\eqref{eq:ODE_SIR} with parameter choices given in Table~\ref{tab:parameters} is shown in Figure~\ref{fig:sample-ODE-and-ABM-run}~(a). For an introduction to mathematical methods for this kind of models we refer the reader to~\cite{DeuflhardBornemann2002}.

\begin{remark} \label{rem:infection-rate}
    For the case of a varying population of size $N$, the infection rate $r_{*\to \circ}$ is given by $r_{*\to \circ} = \bar r_{\mathrm{*\circ}} / N$, where $\bar r_{\mathrm{*\circ}}$ denotes the (population-independent) infection rate. In the literature, $r_{*\to \circ}$ is occasionally also referred to as \emph{transmission rate} or \emph{transmission coefficient}, while $\bar r_{\mathrm{*\circ}}$ is referred to as \emph{infection rate}. For the sake of better readability, we continue to use the term \emph{infection rate} for both rates whenever it is clear from the context.
\end{remark}

\begin{remark}
    There exist a variety of different deterministic models based on differential equations for describing epidemic spreading, for example, ODEs extending the basic SIR model by additional compartments such as hospitalization~\cite{Keeling2021,Wulkow21} or models based on delay differential equations~\cite{DDE20}.
\end{remark}

\subsection{Stochastic Agent-Based Models}
\label{sec:ABM}

The numbers of susceptible, infected, and recovered individuals and their evolution over time can also be computed using an ABM. Often, individuals in the population not only hold an infection status, but are explicitly represented and characterized by, for example, their daily activities in certain locations, as well as other relevant characteristics such as their age or sex.

In simulations, infection events occur with given probabilities whenever an infected and (at least one) susceptible agent meet. The so-called \emph{contact network} arises from the places where agents meet, such as work places, schools, or public transport. For each simulation, at each time step $t$, the numbers $S(t)$, $I(t)$, and $R(t)$ are obtained by counting agents in each infection status. Given the probabilistic nature of infection events (and possibly further elements of the ABM, e.g., agents randomly departing from their activity schedules), these numbers are random variables. Usually, many simulations are performed in order to compute statistical moments, in particular mean and standard deviation, which are taken into account when forecasting the spread of an infection.

\subsubsection{Example 1 -- The Georeferenced Demographic Agent-Based Model}
\label{sec:GERDA-abm}

As our first example, we consider the individualized \textbf{GE}o\-\textbf{R}ef\-er\-enced \textbf{D}e\-mo\-graph\-ic \textbf{A}gent-based mod\-el (GERDA) for the transmission of SARS-CoV-2 and the disease dynamics of COVID-19~\cite{goldenbogen+al2020}. This model uses detailed location data (including homes, workplaces, schools and public places) and a synthetic population of agents with realistic daily schedules to simulate contacts between people in the given locations and the resulting probabilistic infection events. From this point of view, it can be called a \emph{heterogeneous} ABM. The disease progression of infected agents is modeled in detail with rates and times, e.g., for transition to diagnosed, hospitalized, ICU, and recovered or deceased states, depending on the agents' age group. The model has been used with data for several towns (e.g., Tepoztlán (Mexico), Zikhron Ya'akov (Israel), Heinsberg (Germany)). In this work we use the data calibrated for the German municipality of Gangelt, where a well-studied COVID-19 outbreak occurred in early 2020. Based on this model, different types of phar\-ma\-ceu\-ti\-cal and non-phar\-ma\-ceu\-ti\-cal interventions have been analyzed, such as different scenarios of lockdown and reopening, or different vaccination rates and strategies in a population. For further details see~\cite{goldenbogen+al2020}.

We simplify the observation of the disease dynamics computed by GERDA to the basic SIR structure described in Section~\ref{sec:basicODE} by aggregating the different subclasses of infected agents, assuming that all infected agents recover. Additionally, agents are grouped by age, i.e., children (age 0 to 18 years) and adults (19\texttt{+}). Figure~\ref{fig:sample-ODE-and-ABM-run}~(b) shows such an aggregated trajectory in terms of numbers of susceptible, infected and recovered adults and children for a single simulation of GERDA.

\subsubsection{Example 2 -- A Homogeneous ABM}
\label{sec:mock-abm}

In addition to GERDA, we present another ABM for modeling the spread of an infectious disease, but at a coarser level of detail. We assume two homogeneous groups of agents, i.e., adults and children, and that the contact network is complete, i.e., any agent can be infected by any other agent at any time. Thus, we call it a \emph{homogeneous} ABM and use the shorthand notation H/ABM to refer to it. Homogeneity allows tracking only the number of agents per compartment, making the evaluation several orders of magnitude faster than for GERDA. We represent the H/ABM as a Markov jump process that builds on the identical first- and second-order reactions as the ODE model~\eqref{eq:ODE_SIR}. We wish to make explicit the character of an ABM and therefore write the infection rates $r_{*\to \circ}$, for $*,\circ\in\{\mathrm{a},\mathrm{c}\}$, as in Remark~\ref{rem:infection-rate} for a population of size $N$. Additionally, since in GERDA agents can lose their immunity, we introduce a third transition $R_\circ \rightarrow S_\circ$, i.e., recovered individuals from group $R_\circ$ move to the compartment $S_\circ$ and become susceptible again. We assume that this happens independently of the age group with rate $\mu r_\circ$, for $\mu \in \mathopen] 0, 1 \mathclose[$. Thus, the H/ABM represents an SIRS model. Note that in this work the ODE model and the H/ABM share the same parameters. 

The H/ABM can be simulated using Gillespie's stochastic simulation algorithm~\cite{Gillespie1976}, which constructs exact realizations of the H/ABM in continuous time. Furthermore, assuming convergence of $ r_{\mathrm{*} \to \mathrm{\circ}}$ for $N \to \infty$, it can be be shown that the $N$-scaled Markov jump process given by the H/ABM converges to
a frequency-based process, which can be given either by an ordinary or by a stochastic differential equation, see~\cite{Kurtz1978} for details. For example, setting $\mu = 0$, the H/ABM converges for $N \to \infty$ to an ODE of type~\eqref{eq:ODE_SIR}. 

\subsection{Model Agreement} 
\label{sec:agreement}

The ODE model and ABM are supposed to describe the same phenomenon, though at different levels of detail. While many parameters in an ABM can, at least in principle, be measured or observed, e.g., frequencies of people meeting, the six ODE parameters $r_{* \to \circ}$ and $r_*$ (again, for $*,\circ\in\{\mathrm{a},\mathrm{c}\}$) are not directly related to measurable physical quantities. A reasonable parameterization is necessarily based on the agreement of ODE predictions with reality, where ``reality'' can refer to real-world data of infected people or to predictions of a different model taken as a reference (in this work GERDA). We are interested in the latter case and therefore fit the parameters $r_{*\to \circ}$ and $r_*$ in the sense of $L^2$ similarity between the ODE and GERDA results for a fixed number of agents $N_0$.

\begin{figure}[t]
    \centering
    \includegraphics[scale=1]{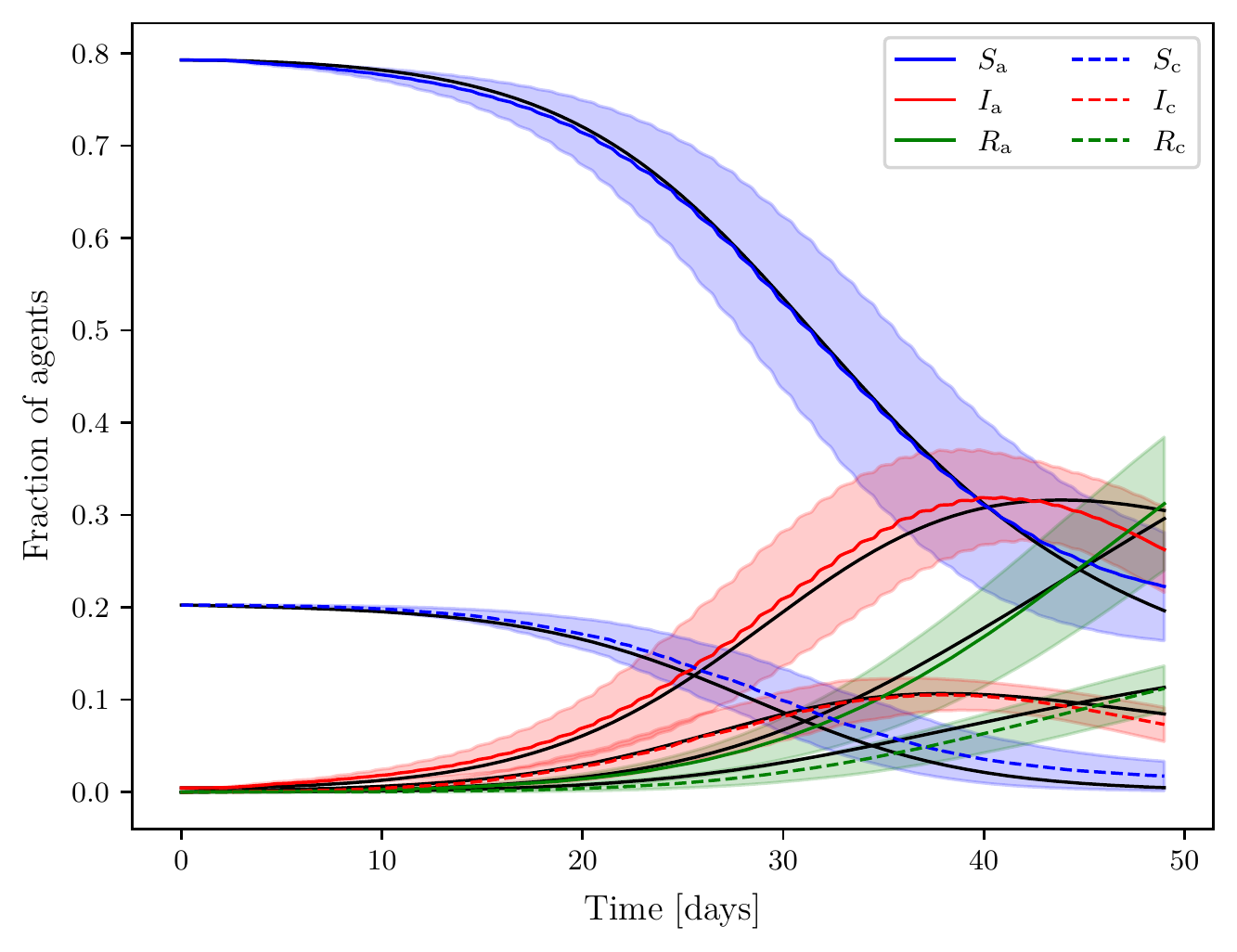}
    \caption{Solution of the fitted ODE system~\eqref{eq:ODE_SIR} (black, solid lines) compared to expected aggregated trajectory of GERDA (mean over 1\,000 independent simulations, sample standard deviation indicated as shaded areas) for parameters given in Table~\ref{tab:parameters}.}
    \label{fig:ABM_ODE_fit}
\end{figure}

A least-squares approach leads to the optimization problem
\begin{equation} \label{eq:least-squares}
    \min_{r_{*\to \circ}, r_*} \int_0^T \left( V_{I,\mathrm{a}}^{-1} (I_\mathrm{a}^{\rm ODE}-I_\mathrm{a}^{\rm ABM})^2 
                                           + V_{I,\mathrm{c}}^{-1} (I_\mathrm{c}^{\rm ODE}-I_\mathrm{c}^{\rm ABM})^2 \right)\, \mathrm{d}t,
\end{equation}
where $I_*^{\rm ABM}$ is computed as a mean of a set of GERDA simulations due to its stochastic nature, and the sample variances $V_{I,*}$ of the respective infected age groups serve as weighting factors for the fit. In Figure~\ref{fig:ABM_ODE_fit} we show the result, obtained by using a nonlinear least-squares solver (an in-built MATLAB solver) for problem~\eqref{eq:least-squares}. The GERDA data is estimated from 1\,000 independent simulations given the parameters in Table~\ref{tab:parameters}. For comparability of GERDA and the H/ABM, we parameterize the H/ABM using the parameters obtained from~\eqref{eq:least-squares}.

\begin{remark} \label{rem:time-resolution}
    The attentive reader will notice in Figure~\ref{fig:ABM_ODE_fit} the undulating trajectories of the GERDA model for susceptible and infected agents, especially in the adult age group. This is due to the fact that the timescale in the GERDA model is measured in hours with different amounts of infection events during night and day.
\end{remark}

%==============================================================================================

\section{Policy Optimization} 
\label{sec:optimization}

When using epidemiological models not only to forecast the spread of infections but also to design favorable policies, three elements are needed in addition to a reliable model: 1)~a design space of policies affecting the model dynamics, e.g., by changing parameters, 2)~an objective to be minimized, and 3)~optimization algorithms. In the following, we will discuss these aspects on the example of the models introduced above.

\subsection{Policies}
\label{sec:policies}

It is well known that the spread of infectious diseases can be influenced not only by phar\-ma\-ceu\-ti\-cal interventions, such as medical treatment or vaccination, but also by non-phar\-ma\-ceu\-ti\-cal interventions, such as social distancing or lockdowns to reduce physical contact between individuals. For our discussion, we consider two time-dependent controls $u_\mathrm{s}, u_\mathrm{w} : [0,T] \rightarrow [0,1]$ as examples representing policies for closing schools and workplaces, or home office for those jobs where work can be done at home, respectively. More precisely, the value of $u_\mathrm{s}(t)$ refers to the fraction of children not going to school at time $t$, i.e., $u_\mathrm{s}(t) = 0$ stands for the policy that no schools are closed and $u_\mathrm{s}(t) = 1$ indicates that all schools are closed. Analogously, the control $u_\mathrm{w}(t)$ is the fraction of workplaces closed. We denote $u = [u_\mathrm{s}, u_\mathrm{w}]^\top$ to refer to the combination of controls. The controls are assumed to be piece-wise constant on a time grid $0 = t_0 < \ldots < t_m=T$ to reflect  the fact that policy changes take some time to be implemented and that continuously changing policies would be complicated to communicate and follow. Consequently, the admissible set of policies
\[
U=\{u:[0,T]\to [0,1]^2 \mid \forall i\in\{1,\dots,m\}: u|_{\mathopen]t_{i-1},t_i\mathclose[} = \mathrm{const} \},
\]
of dimension $n_\mathrm{u} = 2m$ is isomorphic to $[0,1]^{2m} \subset \R^{n_\mathrm{u}}$ such that any control can be represented by $n_\mathrm{u}$ scalar control values. Policies affect the simulated spread of infection by changing the ODE (H/ABM) parameters or, in case of GERDA, by changing the agents' schedules (e.g., staying home instead of going to school). Thus, identical policies must be implemented in different ways depending on the model used.

\subsubsection{Implementing Policies in the ODE}
\label{sec:policy-ode}

In the ODE model presented in Section~\ref{sec:basicODE}, infection rates are determined by the properties of the virus and the population. Additionally, we assume that they depend on the controls $u_\mathrm{s}$ and $u_\mathrm{w}$. For the infection rates within one age group $r_{* \to *}$, for $* \in \{\mathrm{a}, \mathrm{c}\}$, we assume quadratic dependence on the fraction of individuals going to work or school since an infection event can only occur when two individuals meet, i.e.,
\[
    r_{\mathrm{a} \to \mathrm{a}} = r_{\mathrm{aa}} (1-u_{\mathrm{w}})^2, \quad r_{\mathrm{c} \to \mathrm{c}} = r_{\mathrm{cc}} (1-u_\mathrm{s})^2.
\]
Here, the parameters $r_{\mathrm{aa}}$ and $r_{\mathrm{cc}}$ are independent of the policies. Quadratic dependency is widely used for general second-order interactions~\cite{Atkins2008}, but is obviously highly simplistic for epidemiological modeling as, for example, it neglects possible infections outside schools and workplaces. Similarly, we assume the infection rates between age groups $r_{*\to\circ}$ to depend on the controls as well, though with impact reduced by a factor of one half since many such infections will only happen within families, i.e., 
\[
    r_{\mathrm{c} \to \mathrm{a}} = r_{\mathrm{a} \to \mathrm{c}} = r_{\mathrm{ac}} (1-u_{\mathrm{w}}/2) (1-u_{\mathrm{s}}/2).
\]
The recovery rates $r_*$ are assumed to be independent of the controls.

\begin{remark}
    Due to its strong similarity in terms of modeling, the policies for the H/ABM are implemented analogously to the ODE model. For example, using $\Tilde r_{*\circ} \coloneqq r_{*\circ} N_0$, the propensity function that characterizes the infection of an adult by an adult including the policy is given by
    \begin{equation*}
        \tilde r_{\mathrm{a} \to \mathrm{a}} = \frac{\Tilde r_{\mathrm{aa}}}{N} (1-u_{\mathrm{w}})^2,
    \end{equation*}
    where $N$ is the population size used for the H/ABM. The propensity functions $\tilde r_{\mathrm{c} \to \mathrm{c}}$ and $\tilde r_{\mathrm{c} \to \mathrm{a}} = \tilde r_{\mathrm{a} \to \mathrm{c}}$ follow equivalently.
\end{remark}

\subsubsection{Implementing Policies in GERDA}

In contrast to the ODE model, in the GERDA model, where contacts between agents are explicitly modeled, policies do not affect the probability of infection, but the frequency of those contacts. Temporary closure of the workplace or order to home office is done by randomly selecting adults in a ratio equal to the control value $u_{\mathrm{w}}(t)$ and changing their schedules so that they stay home. Analogously, at the rate $u_{\mathrm{s}}(t)$, children stay at home instead of going to school. We assume that an adult from the same household must also be home to supervise children under the age of 13. With a high share of school closures but open workplaces, this leads to an implicit home office obligation for the supervisors involved (i.e., usually the parents). This results in a significant reduction in interactions between agents, which in turn leads to fewer infections and a slower spread of the disease. However, this effect is not reflected by the controls, since they represent only government-mandated interventions.

\begin{remark}
    As a reminder, the intervention polices discussed in this paper are for illustrative purposes only. Further or other policies could be taken, for example, a ban on visiting public places, wearing a mask, contact tracing and isolation, or vaccination.
\end{remark}

\subsection{Objective} 
\label{sec:objective}

For the models and policies introduced in the previous sections, we define a simple example of an objective that nevertheless raises interesting optimization questions. We assume that the general aim of pandemic policy management is to reduce the number $I$ of infected people while minimizing the economic and social costs of the interventions. Let us break down this assumption.

The first aim is reflected in the health objective
\begin{equation*} 
    c_\mathrm{h}(u) = \int_{t=0}^ T (I(t)/N + \exp(10(I(t) - I_\mathrm{max})/N)) \, \mathrm{d}t,
\end{equation*}
where the number of infected agents $I(t)$ depends implicitly on the control $u$ via the ABM or ODE model, respectively. The first term represents the negative impact of infections, which grows approximately linearly with the number of infected people, assuming, for example, that the number of people requiring medical treatment is a fixed fraction of the infected. The second term represents the social impact, which grows drastically as soon as the health care system's capacity to treat those infected who become seriously ill is exhausted. This is formulated in terms of a threshold of infected $I_\mathrm{max}$, again assuming that a fixed fraction of infected agents become seriously ill. 

Both controls incur direct costs, interpreted as economic impact $c_\mathrm{w}$ for home office and social impact $c_\mathrm{s}$ for school closure. The economic impact becomes $+\infty$ when the home office rate approaches an upper bound $u_\mathrm{w}^\mathrm{max} < 1$, reflecting the fact that not all jobs can be stopped or done from home. Thus, we define the second aim $c_\mathrm{w}$ by
\begin{equation*}
    c_\mathrm{w}(u)=  \int_0^T -\log(u_\mathrm{w}^\mathrm{max} - u_\mathrm{w}(t)) \, \mathrm{d}t.
\end{equation*}
We further assume that the social impact of school closures depends quadratically on the fraction of schools closed, i.e.,
\begin{equation*}
    c_\mathrm{s}(u) = \int_0^T  u_\mathrm{s}(t)^2 \,\mathrm{d}t.
\end{equation*}
Combining these three aims with weights $a_\mathrm{s}, a_\mathrm{w} >0$, we define our objective
\begin{equation} \label{eq:objective-sample}
    J(u) \coloneqq c_\mathrm{h}(u) + a_\mathrm{s} c_\mathrm{s}(u) + a_\mathrm{w} c_\mathrm{w}(u).
\end{equation}
Note that a common factor for all three terms does not affect the minimizer, so that we can normalize the weight for $c_h$ to one.

If the number of infected is computed using any time stepping method for solving the SIR initial value problem~\eqref{eq:ODE_SIR}, see, e.g.,~\cite{DeuflhardBornemann2002}, the objective~\eqref{eq:objective-sample} can be evaluated directly using numerical integration. We refer to this as $J^{\rm ODE}$. For an ABM involving randomness, the number of infected is a stochastic process. Each simulation yields a different trajectory for the number of infected. In this case, we refer to~\eqref{eq:objective-sample} as $J^{\rm ABM}$, and define the actual objective as the expectation 
\begin{equation} \label{eq:objective}
    \EE[J^{\rm ABM}(u)].
\end{equation}
To evaluate~\eqref{eq:objective}, a simple Monte Carlo integration 
\begin{equation} \label{eq:objective-MC}
    \EE[J^{\rm ABM}(u)] \approx \frac{1}{n} \sum_{i=1}^n J^{\rm ABM}(u)_i \eqqcolon \EE_n[J^{\rm ABM}(u)]
\end{equation}
based on samples $J^{\rm ABM}(u)_i$ can be used. The accuracy of this unbiased estimate is given in terms of its standard deviation as
\begin{equation} \label{eq:objective-variance}
    \sigma( \EE_n[J^{\rm ABM}(u)]) \approx \frac{\sigma_n(J^{\rm ABM}(u))}{\sqrt{n}}
\end{equation}
with the sample variance 
\[
    \sigma_n^2(J^{\rm ABM}(u)) \approx \frac{1}{n-1} \sum_{i=1}^n (J^{\rm ABM}(u)_i -  \EE_n[J^{\rm ABM}(u)])^2.
\]
The well-known slow convergence of Monte Carlo methods of order $\mathcal{O}(n^{-1/2})$ requires a large number of expensive ABM simulations for a faithful evaluation of the objective. Variance reduction methods such as control variates, quasi-Monte Carlo sequences, or hierarchical Monte Carlo methods are popular approaches to reduce the numerator in~\eqref{eq:objective-variance} and thus the sampling error~\cite{BotevRidder2017}. Unfortunately, they are intrusive by requiring a joint probability space underlying the different random variables used, and thus a carefully designed use of random number generators for sampling. This is difficult to achieve for complex ABMs designed from scratch, and virtually impossible for existing ABMs. Thus, in general, one has to rely on the simple Monte Carlo evaluation~\eqref{eq:objective-MC}.

\subsection{ODE Optimization}
\label{sec:ode-optimization}

For an ODE model, the policy optimization problem
\begin{equation*}
    \min_{u \in U} J^{\rm ODE}(u)
\end{equation*}
is a classical optimal control problem, for which several well-studied numerical solution approaches exist, see, e.g.,~\cite{Betts2010,Rao2009}. We will briefly sketch a simple and straightforward, but for sure not the most efficient, minimization algorithm, which is summarized in Algorithm~\ref{alg:gradient}.

We apply a projected gradient method with line search, where $U$ is the set of admissible controls, i.e., in our case the unit cube in $\R^{n_\mathrm{u}}$. The idea is to go downhill in direction of the steepest descent, which is the negative gradient direction $-\nabla J^{\rm ODE}(u)$ in line~3. This may, however, point to the outside of $U$ if $u_k \in \partial U$. Outward pointing components, for $ i \in \{1, \ldots, n_\mathrm{u}\} $,  are set to zero by the projection 
\begin{equation*}
    P_U(x,y)_i = \begin{cases}
    0, & (y_i=0 \wedge x_i <0) \vee (y_i=1 \wedge x_i>0) \\
    x_i, & \text{otherwise}
    \end{cases}
\end{equation*}
to the admissible set $U$. This search direction $s_k$ is followed downhill for some step size $\alpha$, conceptually the largest step size that 1)~satisfies a constant bound $\alpha \le \alpha_0$ and admissibility $\alpha \le \bar\alpha(u_k) \coloneqq \max\{a \mid u_k+a s_k\in U\}$ in line~4, and 2)~satisfies the Armijo rule of sufficient decrease, 
\begin{equation} \label{eq:armijo}
    J^{\rm ODE}(u_k+\alpha s_k) \le J^{\rm ODE}(u_k) + \alpha c_1 \nabla J^{\rm ODE} (u_k) s_k    
\end{equation}
in line 5. The latter ensures that at least a certain fraction $c_1\in\mathopen]0,1\mathclose[$ of the descent promised by the first-order Taylor approximation of the objective is realized. Algorithm~\ref{alg:gradient} uses a simple backtracking line search in line~6, i.e., the step size for movement in the given direction starts large and is decreased (here halved) until the observed decrease in the objective function is considered sufficient. It therefore can select step sizes that can be up to a factor of two smaller than conceptually desired. The most efficient way to compute the gradient $\nabla J^{\rm ODE}(u)$ in line~3 is to solve the adjoint equations associated with the ODE system~\eqref{eq:ODE_SIR} and $J^{\rm ODE}$. Writing~\eqref{eq:ODE_SIR} compactly as
\[
    \dot y = f(y,u), \quad y = [S_\mathrm{a}, S_\mathrm{c}, I_\mathrm{a}, I_\mathrm{c}]^\top,
\]
where we may safely neglect the recovered compartments, the adjoint equation is the terminal value problem
\[
    -\dot \lambda = f_y (y,u)^\top \lambda + J^{\rm ODE}_y(y)^\top, \quad \lambda(T) = 0,
\]
where the simulated trajectory $y$, i.e., the numerically approximated ODE solution, enters as data.  Then, the gradient is given by 
\[
    \nabla J^{\rm ODE}(u) = \int_0^T f_u(y,u)^\top \lambda \, \mathrm{d}t.
\]
The adjoint equation can be derived from the chain rule using a trivial but clever rearrangement of terms.
For more details we refer to~\cite{GriewankWalther2008}. 

Algorithm~\ref{alg:gradient} converges, usually at a linear rate, to a stationary point, usually a local minimizer, under rather mild regularity assumptions of $J^{\rm ODE}$. Its asymptotic convergence rate deteriorates with growing condition number of the Hessian of $J^{\rm ODE}$, i.e., the ratio of largest and smallest eigenvalue. Preconditioning, i.e., replacing $\nabla J^{\rm ODE}$ by $B^{-1} \nabla J^{\rm ODE}$ with some approximation $B\approx (J^{\rm ODE})''$ of the Hessian, can improve the convergence speed significantly, but constructing effective and computationally cheap preconditioners $B^{-1}$ is not trivial.
For details and more sophisticated algorithms we refer to the nonlinear optimization literature~\cite{NocedalWright2006}.

\begin{algorithm}[t]
    \textbf{Input} $u_0\in U$, $\alpha_0>0$, $c_1 \in \mathopen]0,1\mathclose[$  \\
    \For{$k = 0, \dots, K$}{
        compute $s_k = P_U(-\nabla J^{\rm ODE}(u_k),u_k)$ \\
        $\alpha = \min\{\alpha_0, \bar\alpha(u_k) \}$ \\
        \While{Armijo condition~\eqref{eq:armijo} violated}{
            $\alpha = \alpha / 2$
        }
        $u_{k+1} = u_k + \alpha s_k$
        } 
    \textbf{Output} $u^\star = u_K$
    \caption{Basic steepest descent algorithm}
    \label{alg:gradient}
\end{algorithm}

\subsection{ABM Optimization}
\label{sec:ABM-optimization}

Solving the policy optimization problem for realistic epidemiological ABMs such as the ones presented in Section~\ref{sec:ABM}, i.e.,
\[
    \min_{u\in U} \EE[J^{\rm ABM}(u)],
\]
is much more difficult than solving it for the ODE model due to three major challenges:
\begin{enumerate}[itemsep=0ex, topsep=0.5ex, label=(\roman*)]
    \item ABMs are much more complex, involving thousands of agents, and they often simulate the agents' activities with a higher time resolution, whereas ODE models often only represent the average population and thus omit the vast majority of details that ABMs provide (see Remark~\ref{rem:time-resolution}). As a result, computing a single ABM trajectory is orders of magnitude more expensive.
    \item ABMs are inherently stochastic, such that many independent simulations are required to approximate the objective with sufficient accuracy.
    \item ABMs are inherently discontinuous due to the discrete decisions of the agents. Thus, single ABM trajectories are not differentiable with respect to the control $u$, so that efficient adjoint gradient computation, as possible for ODE models, is not directly available.
\end{enumerate}
Therefore, for the time being, ABM optimization must rely on objective samples alone, inferring descent directions from noisy objective evaluations. We briefly discuss gradient approximation in Section~\ref{sec:gradient-evaluation}, followed by inexact gradient methods in Section~\ref{sec:inexact-gradient-descent}.

\subsubsection{Gradient Approximation}
\label{sec:gradient-evaluation}

For approximating directional derivatives of the objective $\EE[J^{\rm ABM}(u)]$ in direction $v\in\R^{n_\mathrm{u}}$, we use finite differences 
\begin{align}
    \nabla\EE[J^{\rm ABM}(u)]^\top \, v 
    &\approx (2h)^{-1} \left( \EE_n[J^{\rm ABM}(u+hv)]- \EE_n[J^{\rm ABM}(u-hv)]\right) \notag \\
    &= \frac{1}{2hn} \sum_{i=1}^n \left( J^{\rm ABM}(u+hv)_i - J^{\rm ABM}(u-hv)_i \right) \notag \\
    &\eqqcolon \nabla_h  \EE_n[J^{\rm ABM}(u)]^\top \, v \label{eq:finite-differencing}
\end{align}
for small $h>0$. Using the unit vectors $v=e_k$ for $k=1,\dots,n_\mathrm{u}$, the complete gradient vector can be obtained. The approximation error is of order
\begin{equation*}
    \Vert \nabla \EE[J^{\rm ABM}(u)] - \nabla_h  \EE_n[J^{\rm ABM}(u)] \Vert = \mathcal{O}\left(h^2 + \frac{\sigma}{h\sqrt{n}}\right),
\end{equation*}
which requires a careful choice of sample size $n$ of order $n=\mathcal{O}(h^{-6})$ to balance discretization and sampling error. Using 
\[ 
\nabla_{h,k} J^{\rm ABM}(u)_i \coloneqq (2h)^{-1}(J^{\rm ABM}(u+he_k)_i - J^{\rm ABM}(u-he_k)_i),
\]
where $e_k$ is the $k$-th unit vector, in the estimator 
\begin{equation*}
    (V_n)_{kl} = \frac{1}{n-1} \sum_{i=1}^n \left(\nabla_{h,k} J^{\rm ABM}(u)_i - \nabla_h\EE[J^{\rm ABM}(u)]_k\right) \left(\nabla_{h,l} J^{\rm ABM}(u)_i - \nabla_h\EE[J^{\rm ABM}(u)]_l\right)
\end{equation*}
for the sample covariance $V_n \in \R^{n_\mathrm{u}\times n_\mathrm{u}}$, the sample mean standard deviation
\[
    \sigma(\nabla_h J^{\rm ABM}(u)) = \sqrt{\|V_n\| / n},
\]
provides an a posteriori error estimate for the gradient evaluation. Here, $\Vert\cdot\Vert$ refers to the 2-norm. A relative accuracy $\epsilon>0$ of the approximate gradient $\nabla_h E_n[J^{\rm ABM}(u)]$, i.e.,
\begin{equation}\label{eq:relative-gradient-error}
    \| \nabla \EE[J^{\rm ABM}(u)] - \nabla_h \EE_n[J^{\rm ABM}(u)]\| < \epsilon \|\nabla_h \EE_n[J^{\rm ABM}(u)]\|,
\end{equation}
can be ensured with high probability if $n$ is chosen such that $2\sigma \le \epsilon \|\nabla_h \EE_n[J^{\rm ABM}(u)]\|$.
If applicable, variance reduction methods can be used to reduce the sample mean standard deviation $\sigma$ and thus the sampling error. For example, correlated sampling uses the same random number generator seeds for evaluating $J^{\rm ABM}(u-hv)_i$ and $J^{\rm ABM}(u+hv)_i$ in~\eqref{eq:finite-differencing}, which leads to the two values being increasingly correlated for $h\to 0$ and therefore to an expected error order of 
\begin{equation*}
    \Vert \nabla\EE[J^{\rm ABM}(u)] - \nabla_h  \EE_n[J^{\rm ABM}(u)] \Vert = \mathcal{O}\left(h^2 + 1/\sqrt{n}\right).
\end{equation*}
If this increased correlation actually is realized by the ABM, the required number $n$ of samples is reduced significantly to $\mathcal{O}(h^{-4})$. Moreover, if the policy change $v$ affects only policies after some time $t_v>0$, the trajectories corresponding to $u-hv$ and $u+hv$ coincide on $[0,t_v]$ and can be simulated just once using a checkpoint-restart ability of the ABM, if present.

\subsubsection{Inexact Gradient Descent} 
\label{sec:inexact-gradient-descent}

Unlike for the (deterministic) ODE optimization, the exact steepest descent direction is not available when minimizing the stochastic ABM objective~\eqref{eq:objective}. Nevertheless, descent methods still converge as before if sufficient local decrease can be achieved, i.e., if the descent condition $g^\top s \le -c_0 \|g\| \|s\|$ holds for some $c_0>0$, where we write $g \coloneqq \nabla\EE[J^{\rm ABM}(u)]$ for the exact gradient and $s\coloneqq -\nabla_h\EE_n[J^{\rm ABM}(u)] \approx -g$ for the search direction.

A sufficiently small relative error $\epsilon<1/2$ of the computed steepest descent direction, i.e.,  $\|s+g\| \le \epsilon \|s\|$, guarantees 
that $s$ is a descent direction, since by $\|g\| = \|g+s-s\| \ge (1-\epsilon) \|s\|$ it holds that 
\begin{align*}
    g^\top s 
    = g^\top (-g + s+g)   
    \le - \|g\| (1-\epsilon)\|s\| + \epsilon \|g\|\|s\| 
    \le - (1-2\epsilon) \|g\| \|s\|.
\end{align*}
 In particular, choosing $\epsilon = 1/4$ in~\eqref{eq:relative-gradient-error} yields $c_0 = 1/2$ and thus guarantees at least half of the local progress compared to exact steepest descent. Since $-s^\top g \le (1+\epsilon)\|s\|^2$ (using Cauchy--Schwarz) yields
 \[
    -(1+\epsilon) c_1 \alpha \|s\|^2 \le c_1 \alpha s^\top  g,
 \]
 the Armijo condition~\eqref{eq:armijo} is implied by 
 \[
    \EE[J^{\rm ABM}(u+\alpha s)] - \EE[J^{\rm ABM}(u)] \le -(1+\epsilon) c_1 \alpha \|s\|^2.
 \]
 When evaluating the objective with absolute sampling errors $e \le 2\sigma$ as in~\eqref{eq:objective-variance}, we obtain
 \[
  \EE[J^{\rm ABM}(u+\alpha s)] - \EE[J^{\rm ABM}(u)] \le \EE_n[J^{\rm ABM}(u+\alpha s)] - \EE_n[J^{\rm ABM}(u)] + 2e. 
 \]
 Choosing $e\le \epsilon c_1\alpha\|s\|^2$, the Armijo condition is satisfied with high confidence if the acceptance test
 \begin{equation}\label{eq:inexact-acceptance-test}
    \EE_n[J^{\rm ABM}(u+\alpha s)] - \EE_n[J^{\rm ABM}(u)] \le -(1+3\epsilon) c_1 \alpha \|s\|^2
 \end{equation}
 is passed for the trial step $s$.
 
 \begin{algorithm}[t]
    \textbf{Input} $u_0\in U$, $\alpha_0>0$, $c_1\in\mathopen]0,1\mathclose[$ \\
    \For{$k = 0, \dots, K$}{
        compute $s_k = P_U(-\nabla_h\EE_n[J^{\rm ABM}(u_k)],u_k)$ with relative accuracy $\epsilon < 1/2$ \\
        $\alpha = \min\{\alpha_0,\bar\alpha(u_k)\}$ \\
        \While{acceptance test~\eqref{eq:inexact-acceptance-test} fails}{
            $\alpha = \alpha / 2$
        }
        $u_{k+1} = u_k + \alpha s_k$
        } 
        \textbf{Output} $u^\star = u_k$
    \caption{Basic inexact gradient descent algorithm}
    \label{alg:inexact-gradient}
\end{algorithm}

Algorithm~\ref{alg:inexact-gradient} converges similarly to the exact gradient descent to a stationary point, usually a local minimizer. However, it is computationally expensive due to the large number of samples required to compute approximate gradients (line~3) and for the line search (line~5).

\begin{remark}
    An alternative is stochastic gradient descent, where instead of averaging gradient samples as in~\eqref{eq:finite-differencing} and then taking an inexact gradient step, the gradient averaging can be done implicitly on the fly, by taking sufficiently small steps in the direction of a single or an average of a few samples. This is known as finite difference based stochastic approximation or Kiefer--Wolfowitz algorithm~\cite{KieferWolfowitz1952}, a variant of which is given in the supplementary material. Its convergence rate is in general sublinear. Several improvements have been proposed, such as simultaneous perturbation gradient approximation~\cite{Spall1992} reducing the effort of gradient sampling from $\mathcal{O}(n_\mathrm{u})$ to $\mathcal{O}(1)$. For an overview of improved stochastic gradient descent approaches we refer to~\cite{ChauFu2015,BottouCurtisNocedal2018}. However, even with these improvements, the convergence is sublinear, requiring many iterations.
\end{remark}

\subsection{Multilevel Optimization} 
\label{sec:ml-opt}

Nonlinear multilevel or multigrid optimization approaches such as MG/OPT~\cite{Nash2000} or recursive multilevel trust region (RMTR)~\cite{GrattonSartenaerToint2008} exploit the computationally cheap optimization of coarser models to generate trial steps that can make much more progress than possible with Taylor approximations of the fine model by respecting the underlying nonlinear structures. In the case of an ABM as fine model, a well-matched ODE model as given in Section~\ref{sec:basicODE} can serve as a coarse model, with the advantage of being computationally extremely cheap. Algorithm~\ref{alg:RMTR} implements a basic multilevel optimization algorithm. Here, it is restricted to only two levels $J^{\rm ABM}$ and $J^{\rm ODE}$.

\begin{algorithm}[t]
    \textbf{Input} $u_0\in U$, $\rho_0>0$ \\
    \For{$k = 0, \dots, K$}{
        compute $s_k = P_U(-\nabla_h\EE_n[J^{\rm ABM}(u_k)],u_k)$ with relative accuracy $\epsilon < 1/2$ \\
        $\rho = \rho_0$ \\
        \Repeat{acceptance test~\eqref{eq:multilevel-acceptance-test}}{
            $\delta u_k = \mathop\mathrm{arg\, min}\limits_{\delta u \in (U-u_k) \cap \mathopen]-\rho,\rho\mathclose[^{n_\mathrm{u}}} 
                    J^{\rm ODE}(u_k+\delta u) - (s_k + \nabla J^{\rm ODE}(u_k))^\top  \delta u$ \\   
            $\rho = \rho/2$ \\
        }
        $u_{k+1} = u_k + \delta u_k$
    } 
    \textbf{Output} $u^\star = u_k$
    \caption{Basic multilevel optimization}
    \label{alg:RMTR}
\end{algorithm}
The two main differences to the inexact gradient descent algorithm are the use of a coarse model for computing the trial step instead of taking the negative gradient direction, and by restricting the trial step to an $l^\infty$ trust region within the $n_\mathrm{u}$-dimensional box $U$ of admissible policies instead of performing a line search. For minimizing the coarse ODE model in line~6, the deterministic steepest descent algorithm (Algorithm~\ref{alg:gradient}) is a viable option. Compared to evaluating ABM gradients, the solution of an ODE optimization problem incurs negligible computational cost.

In order to guarantee convergence, the step suggested by optimizing the coarse model must lead to a sufficient decrease of the fine objective as required by the Armijo type condition~\eqref{eq:inexact-acceptance-test} for some trust region radius $\rho>0$. More precisely, the acceptance test~\eqref{eq:inexact-acceptance-test}, which reformulates for the multilevel optimization scheme to
\begin{equation}
    \EE_n[J^{\rm ABM}(u + \delta u_k)] - \EE_n[J^{\rm ABM}(u)] \le -(1+3\epsilon) c_1 \delta u_k^\top s \label{eq:multilevel-acceptance-test}
\end{equation}
with $e \le \epsilon c_1 \delta u_k^\top s$, needs to be fulfilled. This can be ensured if the coarse model is first-order consistent, i.e., its gradient coincides with the fine model gradient at the current iterate $u_k$, which is achieved by adding the linear correction $-(s_k + \nabla J^{\rm ODE}(u_k))^\top  \delta u$ in line~6. Even then, the computed step does not have to lead to a decrease of the fine objective $J^{\rm ABM}$ since the coarse and fine model are different. A backtracking line search is not guaranteed to solve this issue since in highly nonlinear problems the step $\delta u_k$ does not necessarily have to be a descent direction. Instead, the minimization of the coarse model is restricted to a neighborhood of the current iterate $u_k$. Using the $l^\infty$ neighborhood $u_k + \mathopen]-\rho,\rho\mathclose[^{n_\mathrm{u}}$ and intersecting it with the admissible region $U$ leads to the same type of box-constrained optimization subproblem as before. For sufficiently small $\rho>0$, a reduction of $J^{\rm ABM}$ is ensured by first-order consistency. Note that the selection of the trust region radius $\rho$ in Algorithm~\ref{alg:RMTR} is deliberately simple. For more efficient update strategies we refer to~\cite{ConGouToi2000}.

Also note that every iteration of the multilevel algorithm is about as expensive as one of the inexact gradient descent algorithm. An acceleration of optimization by multilevel schemes can only be accomplished by reducing the number of iterations until meeting the desired accuracy. In the early iteration phase, much larger -- and therefore fewer -- steps can be taken if the coarse model captures the global nonlinear structure of the fine objective well, and much better than a first- or second-order Taylor model does. In the later phase, the coarse model's Hessian acts as a preconditioner for the steepest descent method on the fine model. A good agreement of both models' Hessians improves the asymptotic convergence rate on ill-conditioned problems and avoids the ``zig-zagging'' that slows down gradient methods applied to general ill-conditioned problems.

Algorithm~\ref{alg:RMTR} can be extended from two-levels to a true recursive multilevel scheme by solving the coarse model subproblem in line~6 not by a simple gradient method, but by minimizing some even coarser model. However, due to the extreme difference in computational effort between agent-based and ODE models, the use of coarser ODE models does not promise any benefit. In contrast, the use of a hierarchy of fine ABM, coarse ABM, SDE models, and ODE models is an interesting perspective, which, however, we will not explore in this work.

%==============================================================================================

\section{Numerical Examples} 
\label{sec:num-examples}

In this section, we present the results of using the optimization techniques for ABMs introduced in the previous section for the given example objective. In Section~\ref{sec:parameters}, we parameterize the models. We then compare the multilevel optimization~(MLO) algorithm (Algorithm~\ref{alg:RMTR}) with the inexact gradient descent~(IGD) algorithm (Algorithm~\ref{alg:inexact-gradient}), and show the optimization results corresponding to our examples GERDA and the H/ABM in Sections~\ref{sec:optimal-policies_mock} and~\ref{sec:optimal-policies}.

\subsection{Parameterization} 
\label{sec:parameters}

In order to anchor the examples in some empirical orders of magnitude, we considered figures relating to the beginning of the pandemic in Germany, as the model example focuses on the German municipality of Gangelt. In each experiment we consider a time frame of seven weeks, i.e., $T=1\,176$ hours. We consider piecewise constant policies of different durations, i.e., control time grid size. 

\begin{table}[t]
\centering
\caption{Parameters used for the model simulations. Here, the infection rates $r_{\mathrm{*\circ}}$ denote the control independent parameters $r_{\mathrm{*\circ}}$ for the ODE and $\Tilde r_{\mathrm{*\circ}}$ for the H/ABM, respectively. Note that the number of agents for the GERDA model is subject to random changes with each newly generated world.}
\label{tab:parameters}
\begin{tabular}{lccc}
\toprule
                                     & \multicolumn{3}{c}{Model}                                                                                       \\ \cmidrule{2-4} 
Parameter                            & GERDA             & ODE                                          & H/ABM                                        \\ \midrule
Time $T$ {[}hours{]}                 & 1\,176            & 1\,176                                       & 1\,176                                       \\
Number of agents $N$                 & 1\,091            & --                                           & 1\,091                                       \\
Initially infected $I_\mathrm{a}(0)$ & 5                 & 5/1\,091                                     & 5                                            \\
Initially infected $I_\mathrm{c}(0)$ & 0                 & 0                                            & 0                                            \\
World                                & Gangelt (reduced) & --                                           & --                                           \\
General infectivity                  & 0.175             & --                                           & --                                           \\
General interaction frequency        & 1                 & --                                           & --                                           \\
Infection rate $r_{\mathrm{aa}}$     & --                & \multicolumn{1}{l}{$1.0252 \times 10^{-12}$} & \multicolumn{1}{l}{$1.1185 \times 10^{-9}$}  \\
Infection rate $r_{\mathrm{ac}}$     & --                & \multicolumn{1}{l}{$4.8804 \times 10^{-4}$}  & \multicolumn{1}{l}{$5.3246 \times 10^{-1}$}  \\
Infection rate $r_{\mathrm{cc}}$     & --                & \multicolumn{1}{l}{$6.1482 \times 10^{-13}$} & \multicolumn{1}{l}{$6.7077 \times 10^{-10}$} \\
Recovery rate $r_\mathrm{a}$         & --                & \multicolumn{1}{l}{$4.2148 \times 10^{-2}$}  & \multicolumn{1}{l}{$4.2148 \times 10^{-2}$}  \\
Recovery rate $r_\mathrm{c}$         & --                & \multicolumn{1}{l}{$4.3427 \times 10^{-2}$}  & \multicolumn{1}{l}{$4.3427 \times 10^{-2}$}  \\
Immunity loss $\mu$                  & --                & --                                           & 0.2                                          \\ \bottomrule
\end{tabular}
\end{table}

According to~\cite{ritter+al2021}, who analyze intensive care unit loads for Berlin, Madrid and Lombardy, the fraction of infected in Berlin that needed intensive care in the first wave of the pandemic was 6~\%. We use this number as an approximation for Germany and hence Gangelt. With a number of intensive care beds of around 32\,000 in Germany\footnote{From mid April 2020, intensive care beds had to be registered with a central agency. In the second half of April, the numbers of beds registered ranged around this value, see the daily reports for this time span under \cite{rki2021}.}, if 6~\% of the infected should not exceed this number, the maximum of infected that could be allowed would be around 530\,000 persons, which corresponds to about 0.6~\% of the population of about 83 million. As not all intensive care beds can be allocated to COVID-19 patients since there are other reasons for needing such a bed, we reduce this number to 0.5~\%.    
The weights given to the costs of home office and home schooling versus the weight given to the number of infected are subject of a societal or political debate that we do not intend to enter here. As this point is important but beyond the scope of this paper, we reiterate that the chosen objective function serves as an example for the development of optimization methods. Therefore, an example is defined rather than opting for a multicriteria analysis that would leave the choice of weights up to a potential user of the results obtained. Thus, we set $a_\mathrm{s} = a_\mathrm{w}=1$.

For the upper bound on the fraction of work that can be done at home, we choose an educated guess of $u_\mathrm{w}^\text{max} = 0.81$, based on the following considerations. According to~\cite{hammermann+voigtlaender2020}, the share of office jobs is below 36.7~\%, but a lockdown can also affect non-office jobs, e.g., sales or industrial employees, or employees in the service sector. On the other hand, there are jobs that even in a lockdown cannot be closed down or moved to home office. We consider the examples of workers in the health care system (5.7 million employed), food sales (552\,200 employed), as well as police and fire departments (1.7 million, all numbers from \cite{destatis}). These sum to about 8.5 million, representing about 17.5~\% of the approximately 45 million people employed in Germany. To account for further cases not included here, we increase this number to 19~\%, meaning that the economic costs explode when approaching 81~\% of work not done at the workplaces. 

All parameters relevant to the models are summarized in Table~\ref{tab:parameters}, while the parameters relevant to the optimization algorithms are summarized in Table~\ref{tab:opt-parameters}. Every algorithm is started from $u_0 = 0_{n_\mathrm{u}}$, where $0_{n_\mathrm{u}}$ denotes a vector of length $n_\mathrm{u}$ with zeros only.

\begin{table}[t]
\centering
\caption{Parameters used for the optimization algorithms.}
\label{tab:opt-parameters}
\begin{tabular}{lc}
\toprule
Parameter                                              &          \\ \midrule
Fraction of descent $c_1$                              & 0.1      \\
Maximum trust region radius $\rho_0$                   & 0.5      \\
Maximum finite differencing step size $h_\text{max}$   & 0.1      \\
Accuracy $\epsilon$                                    & 0.25     \\
Health care system's capacity threshold $I_\text{max}$ & $0.005N$ \\
Threshold economic impact $u_\mathrm{w}^\text{max}$    & 0.81     \\
Initial $u_0$                                          & $0_{n_\mathrm{u}}$        \\
Weights $a_\mathrm{s} = a_\mathrm{w}$                  & 1        \\ \bottomrule
\end{tabular}
\end{table}

\subsection{Optimal Policies for the H/ABM}
\label{sec:optimal-policies_mock}

Due to the very long runtime of GERDA (a single simulation with the parameters given in Table~\ref{tab:parameters} takes approximately 30 seconds on a laptop), a detailed numerical analysis of the algorithms would not be justifiable for this model. Therefore, in the following, we use the H/ABM introduced in Section~\ref{sec:mock-abm} that has a strongly reduced computational cost to compare the multilevel algorithm with the inexact gradient descent algorithm. As a termination criterion, we choose either 15 iterations or when the sample size $n$ of the objective estimate~\eqref{eq:objective-MC} or gradient estimate~\eqref{eq:finite-differencing} is more than one million samples, whichever comes first.

\begin{figure}[t]
    \centering
    \begin{subfigure}{0.49\textwidth}
        \centering
        \includegraphics[scale=1]{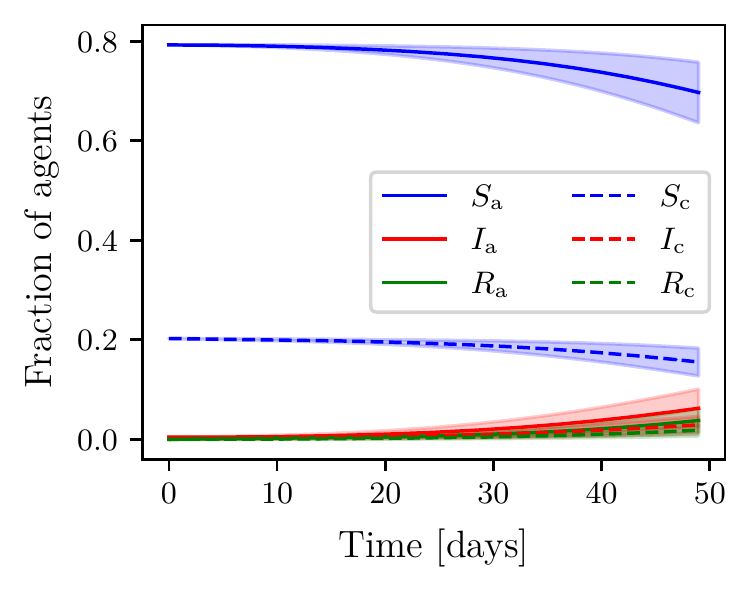}
        \caption{}
    \end{subfigure}
    \begin{subfigure}{0.49\textwidth}
        \centering
        \includegraphics[scale=1]{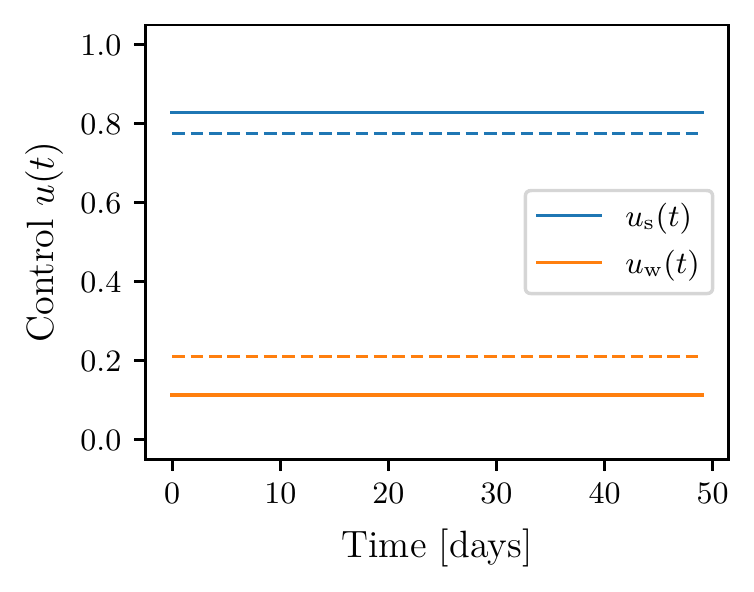}
        \caption{}
    \end{subfigure}
    \begin{subfigure}{0.49\textwidth}
        \centering
        \includegraphics[scale=1]{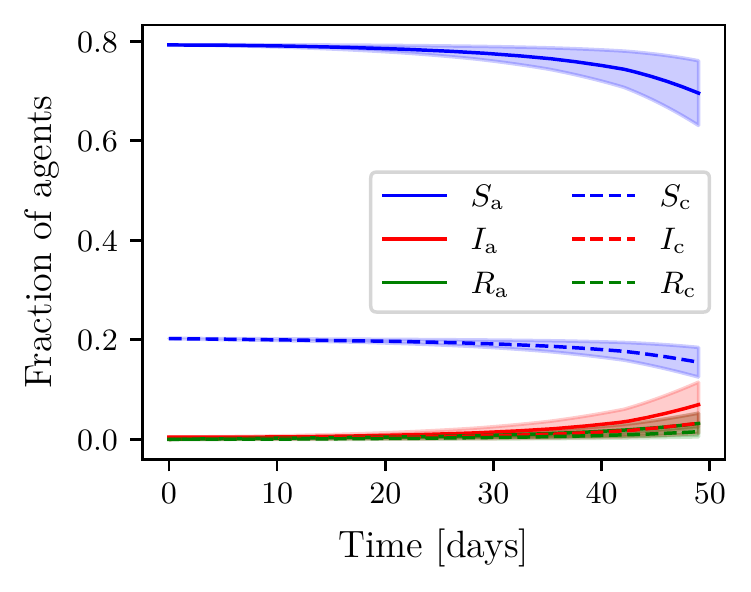}
        \caption{}
    \end{subfigure}
    \begin{subfigure}{0.49\textwidth}
        \centering
        \includegraphics[scale=1]{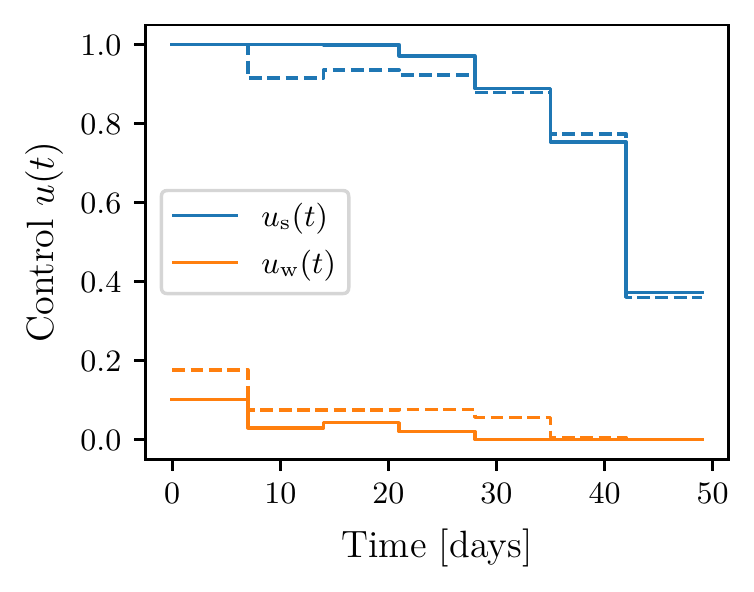}
        \caption{}
    \end{subfigure}
    \caption{Expected time courses of the H/ABM and corresponding optimal policies computed using the multilevel algorithm (solid) compared to the policies computed using the inexact gradient descent algorithm (dashed) for constant and weekly changing policies.}
    \label{fig:dummy_opt_ctrl_cmp}
\end{figure}

First, we present the policies computed by Algorithms~\ref{alg:inexact-gradient} and~\ref{alg:RMTR} for two different cases: constant policies that do not change for the entire time horizon of seven weeks, and piece-wise constant policies that change every 168~hours (i.e., weekly changing). Figure~\ref{fig:dummy_opt_ctrl_cmp} shows the computed policies as well as the corresponding expected time course of the H/ABM for a time frame of 7~weeks. The solid line represents the solution of the multilevel algorithm and the dashed line represents the solution of the inexact gradient descent algorithm. Since both solutions are qualitatively similar, we expect qualitatively similar results in the expected time course for continuity reasons. Both for constant and piece-wise constant policies, the share of home schooling is fairly high, while the share of home office is quite low. For the piece-wise constant policy, the share of home office is close to zero after the first week, but at the expense of home schooling. Independently from the policy chosen, the expected time course of the H/ABM is almost identical, with infections increasing only towards the end of the simulation. 

\begin{figure}[t]
    \centering
    \begin{subfigure}{0.49\textwidth}
        \centering
        \includegraphics[width=\textwidth]{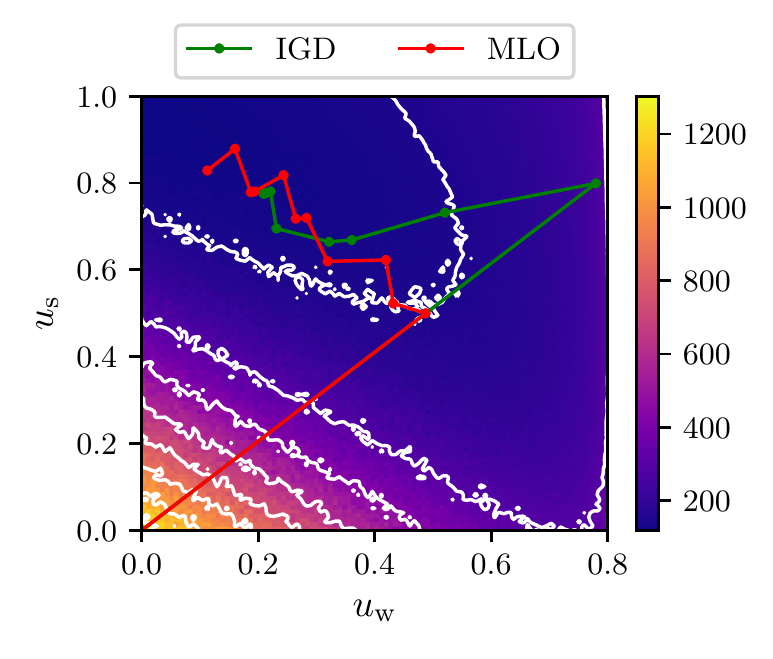}
        \caption{}
    \end{subfigure}
    \begin{subfigure}{0.49\textwidth}
        \centering
        \includegraphics[width=\textwidth]{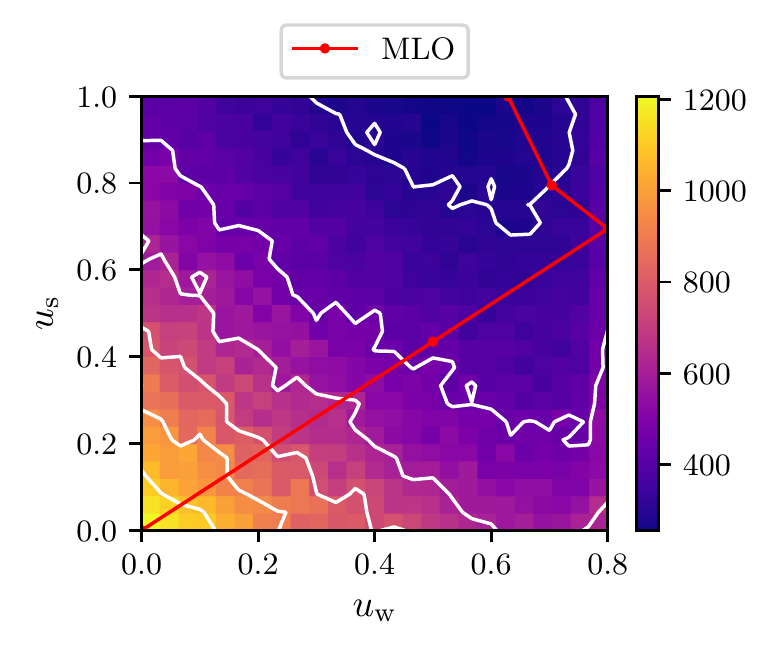}
        \caption{}
    \end{subfigure}
    \caption{Monte Carlo evaluations $\mathbb{E}_n[J^\text{ABM}(u)]$ for $n=100$ and $u=[u_\mathrm{s}, u_\mathrm{w}]^\top \in \mathbb{R}^{n_\mathrm{u}}$ with $n_\mathrm{u}=2$, i.e., the policies are constant for the entire simulation time $T$, and optimization path given by iterates $u_k$ computed using the multilevel and inexact gradient descent algorithms. The objective $J^\text{ABM}$ in Figure~(a) is evaluated for the H/ABM and in Figure~(b) for the GERDA model. For $n_\mathrm{u}=2$ it can be seen clearly that both objectives have different minimizers.}
    \label{fig:energylandscapes}
\end{figure}

\begin{figure}[th!]
    \centering
    \includegraphics[scale=1]{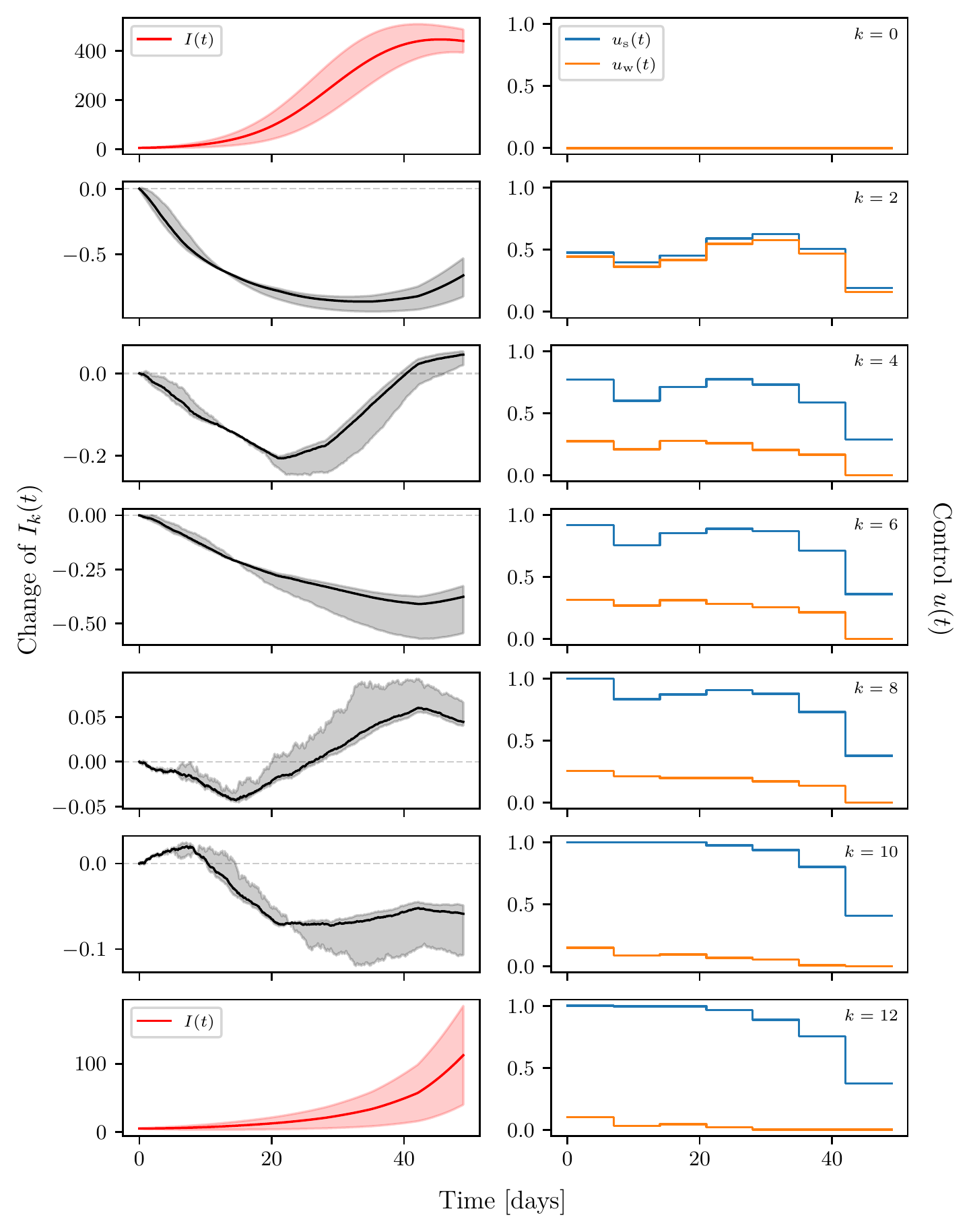}
    \caption{Expected number (red/solid) and standard deviation (red/shaded area) of infected agents $I(t) = I_\mathrm{a}(t) + I_\mathrm{c}(t)$ for the corresponding policies at iteration $k$. Expected relative change (black/solid) and standard deviation (gray/shaded area) to the iteration $k-2$ (i.e., the previous figure). Data estimated from 10\,000 H/ABM simulations.}
    \label{fig:opt-ctrl-progress}
\end{figure}

\begin{figure}[t]
    \centering
    \begin{subfigure}{0.49\textwidth}
        \centering
        \includegraphics[scale=1]{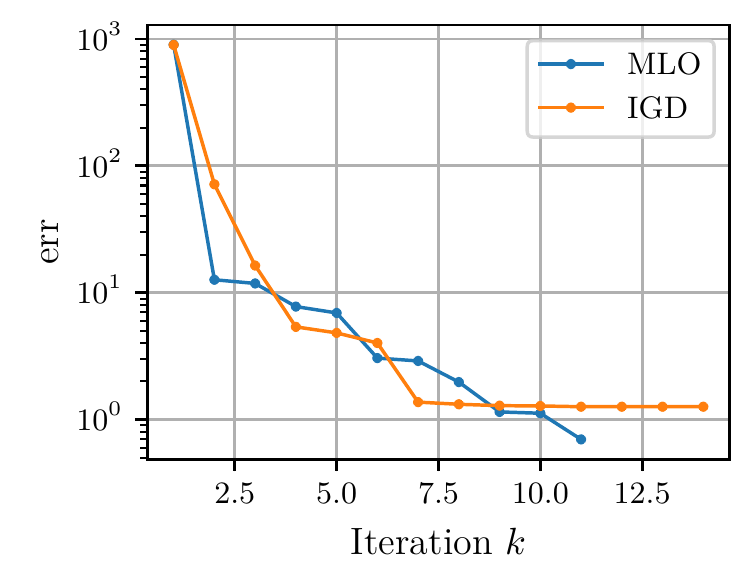}
        \caption{}
    \end{subfigure}
    \hfil
    \begin{subfigure}{0.49\textwidth}
        \centering
        \includegraphics[scale=1]{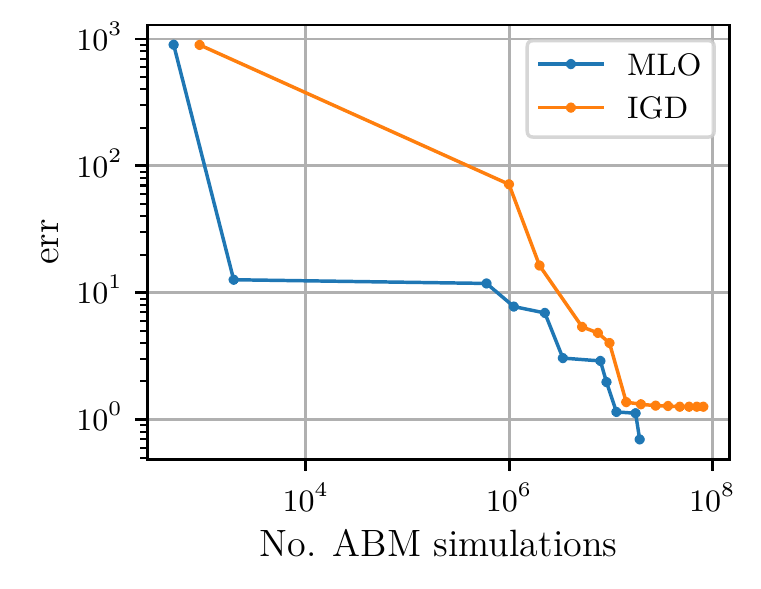}
        \caption{}
    \end{subfigure}
    \begin{subfigure}{0.49\textwidth}
        \centering
        \includegraphics[scale=1]{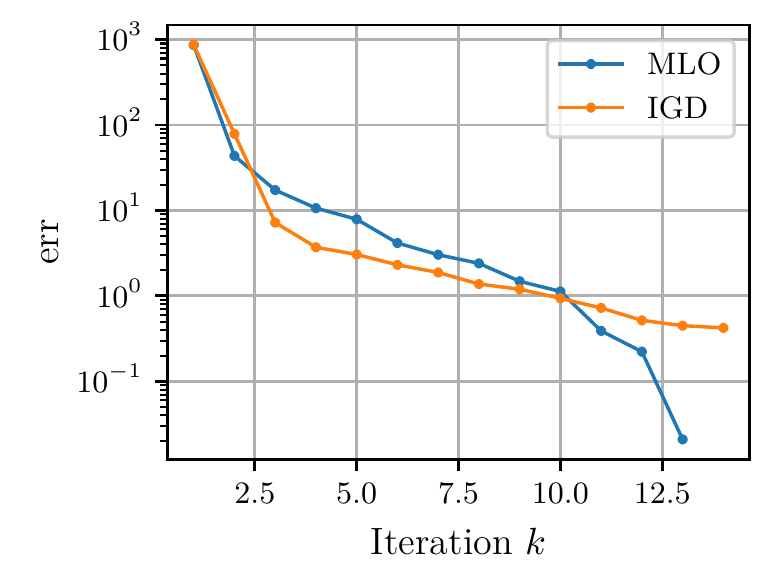}
        \caption{}
    \end{subfigure}
    \hfil
    \begin{subfigure}{0.49\textwidth}
        \centering
        \includegraphics[scale=1]{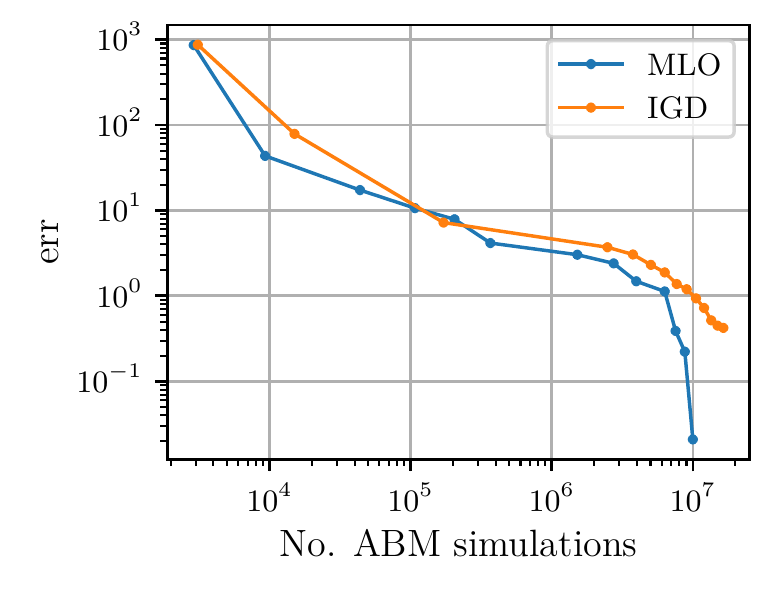}
        \caption{}
    \end{subfigure}
    \caption{Convergence of the multilevel and inexact gradient descent algorithms with respect to the number of iterations and number of ABM simulations for~(a-b) constant and~(c-d) piece-wise constant policies, i.e., weekly changing.}
    \label{fig:convergence-mock}
\end{figure}

\begin{figure}[t]
    \centering
    \begin{subfigure}{0.49\textwidth}
        \centering
        \includegraphics[scale=1]{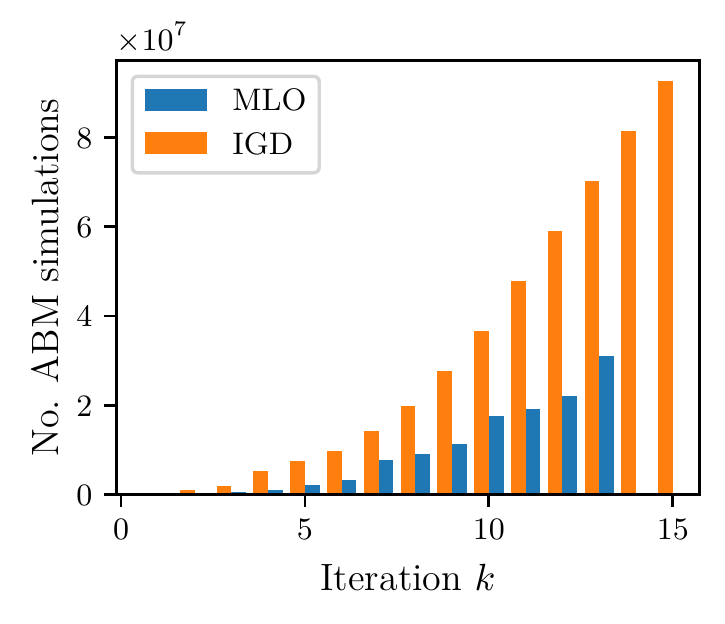}
        \caption{}
    \end{subfigure}
    \begin{subfigure}{0.49\textwidth}
        \centering
        \includegraphics[scale=1]{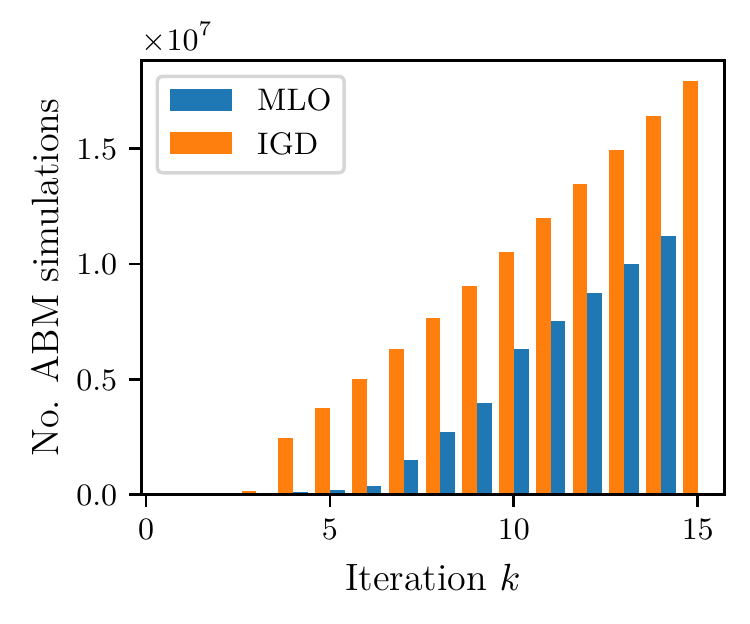}
        \caption{}
    \end{subfigure}
    \caption{Number of ABM evaluations for the H/ABM with~(a) constant controls and~(b) piece-wise constant controls, precisely weekly changing controls. It can be clearly seen that the multilevel algorithm performs better than the inexact gradient descent algorithm at the beginning in terms of the number of iterations. In the further course, both algorithms become expensive at about the same rate.}
    \label{fig:barchart}
\end{figure}

Next, we take a look at the path that both algorithms are taking. Figure~\ref{fig:energylandscapes} shows the Monte Carlo evaluation of $\mathbb{E}_n[J^\text{ABM}(u)]$ for $n=100$ and $u=[u_\mathrm{s}, u_\mathrm{w}]^\top \in \mathbb{R}^2$, i.e., constant policies. Note that an explicit approximation as done here is only possible for $n_\mathrm{u} = 2$. For $n_\mathrm{u} > 2$ it becomes computationally (and visually) infeasible. We see that, as expected, the multilevel algorithm is initially superior to the inexact gradient descent and takes more advantageous steps towards the minimizer $u^\star$ of $J^\text{ABM}$. This is because it respects the underlying nonlinear structures of H/ABM, which the inexact gradient descent algorithm cannot do. The first step of the inexact gradient descent algorithm is just barely still admissible. As it continues, the objectives $J^\text{ABM}$ computed by the algorithms descend, however, at different rates. When approaching the optimum, the changes with respect to $u^\star$ get less pronounced (cf.~Figure~\ref{fig:opt-ctrl-progress}). Especially the inexact gradient descent algorithm slows down significantly and barely makes progress. This can also be seen in Figure~\ref{fig:convergence-mock}, which shows the convergence. We define the error as
\begin{equation*}
    \text{err} \coloneqq \vert \EE[J^\text{ABM}(u_k)] - \EE[J^\text{ABM}(u^\star)] \vert,
\end{equation*}
where $u^\star$ denotes the true minimizer of $J^\text{ABM}$ and $u_k$ the current iterate. The true minimizer $u^\star$ is unknown and can only be approximated. Thus, we set $u^\star = u_{k_\text{max}}^\text{MLO}$, where $u_{k_\text{max}}^\text{MLO}$ denotes the final iterate using the multilevel algorithm algorithm since it can be shown that $J^\text{ABM}$ is locally convex (cf.~Figure~\ref{fig:energylandscapes}~(a)) and that $u_{k_\text{max}} \le u_{k_\text{max}}^\text{IGD}$ holds. The figures confirm that in both cases the multilevel algorithm outperforms the inexact gradient descent algorithm in the early phase. In the later phase the increasing number of ABM simulations is driven by the choice of $e$, which controls the number $n$ required for the estimates of $\EE[J^\text{ABM}(u)]$ to pass the acceptance test with high confidence (cf.~Sections~\ref{sec:inexact-gradient-descent} and~\ref{sec:ml-opt}). Figure~\ref{fig:barchart}, which shows the cumulative sum of ABM simulations, illustrates this very well. It also shows that both algorithms become expensive at about the same rate.

\subsection{Optimal Policies for GERDA}
\label{sec:optimal-policies}

We now consider the GERDA model introduced in Section~\ref{sec:GERDA-abm}. To handle the massive computational effort, we use a reduced version of the municipality of Gangelt with about 1\,000 agents. In this way, one agent in GERDA is roughly equivalent to 10~people in the real world (as of December~2021). This allows computing objectives and gradient estimates as in Section~\ref{sec:ABM-optimization} for large $n$ in reasonable time on a high performance computer.

\begin{figure}[t!]
    \centering
    \begin{subfigure}{0.49\textwidth}
        \centering
        \includegraphics[scale=1]{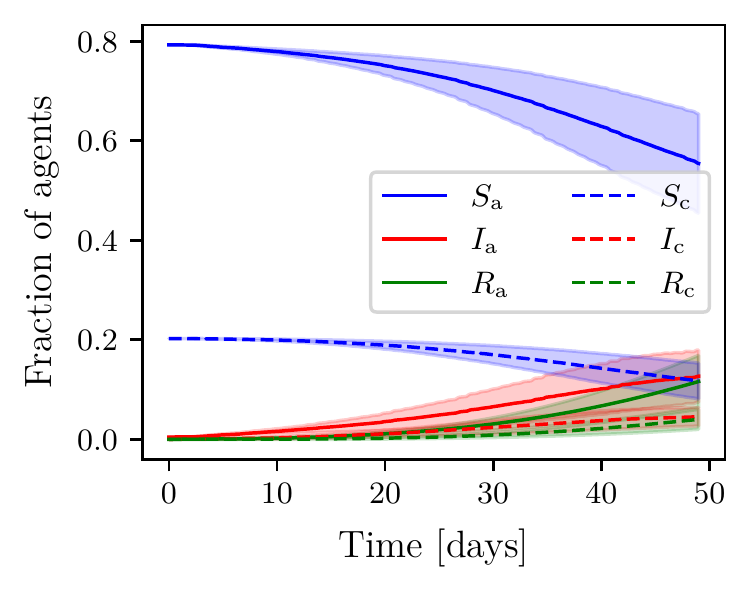}
        \caption{}
    \end{subfigure}
    \begin{subfigure}{0.49\textwidth}
        \centering
        \includegraphics[scale=1]{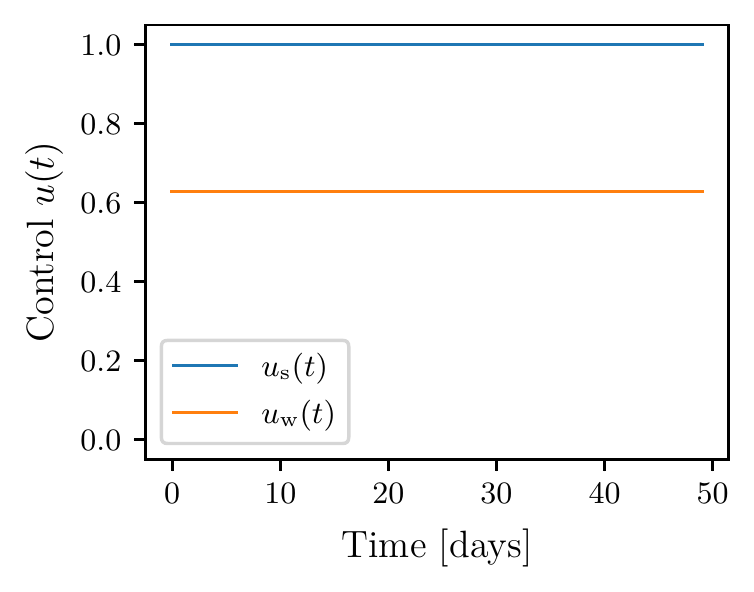}
        \caption{}
    \end{subfigure}
    \begin{subfigure}{0.49\textwidth}
        \centering
        \includegraphics[scale=1]{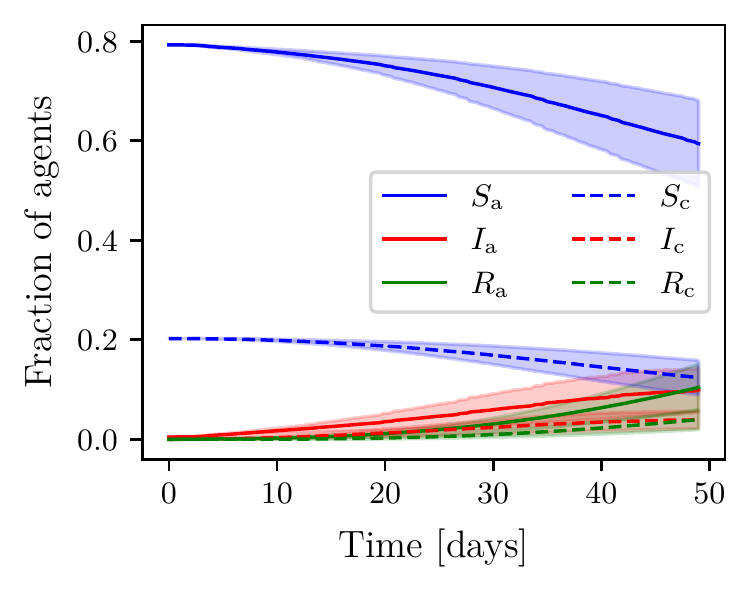}
        \caption{}
    \end{subfigure}
    \begin{subfigure}{0.49\textwidth}
        \centering
        \includegraphics[scale=1]{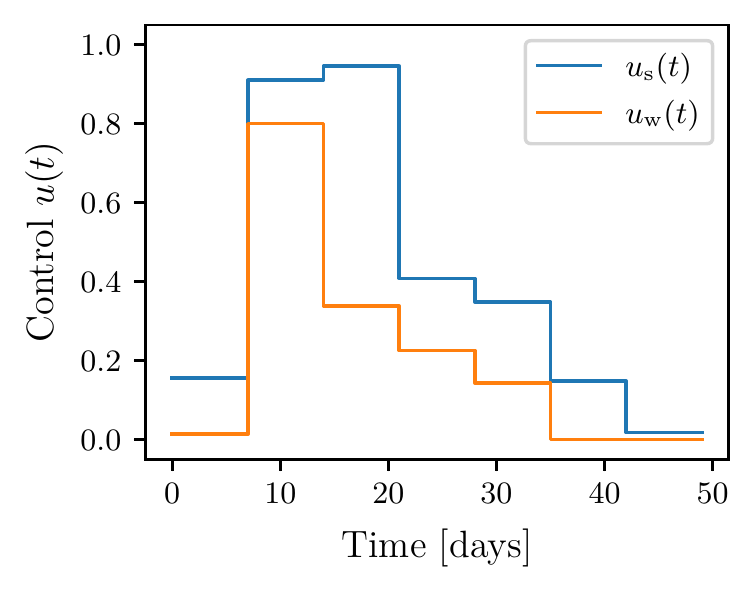}
        \caption{}
    \end{subfigure}
    \caption{Expected time courses of GERDA and corresponding optimal policies computed using the Multilevel Optimization Algorithm~\ref{alg:RMTR} for constant and weekly changing policies. Note that compared to the policies shown in Figure~\ref{fig:dummy_opt_ctrl_cmp} for the H/ABM, the policies for GERDA propose a much stricter lockdown.}
    \label{fig:GERDA_opt_ctrl_cmp}
\end{figure}

Now, we present the results for GERDA computed by the multilevel optimization algorithm for the same two scenarios, i.e., constant and weekly changing policies. Figure~\ref{fig:GERDA_opt_ctrl_cmp} shows the optimal policies for both scenarios. As for the H/ABM, the the difference in the effects of constant and piece-wise constant policies on the trajectories are rather negligible for GERDA. However, despite appropriately adjusted parameterization of the two models (cf.~Section~\ref{sec:agreement}), we notice that the minimizers $u^\star_\text{H/ABM}$ and $u^\star_\text{GERDA}$ are different; compare the solutions shown in Figures~\ref{fig:energylandscapes} or~\ref{fig:GERDA_opt_ctrl_cmp}. More precisely, in both scenarios the policies $u=[u_\mathrm{s}, u_\mathrm{w}]^\top$ are higher, especially the share of home office $u_\mathrm{w}$. Remarkably, for the weekly changing scenario, there is an extreme jump in both policies in the second week of approximately 80 percentage points. The ``delay'' of one week and the subsequent jump are due to the incubation period modeled in GERDA. In the model this results in low values for the policies; in the real world, with infection rates in the low double digits, a hard lockdown in the first seven days as shown in Figure~\ref{fig:dummy_opt_ctrl_cmp} for the H/ABM could lead to resentment or even noncompliance from the public. In the further course, both policies are steadily scaled back. The policies calculated for GERDA are more in line with what was implemented by many governments in early 2020. This shows the need to optimize detailed ABMs directly, since coarser models -- even well-fitted ABMs -- seem to exhibit a crucial loss of information. 

The GERDA model computations emphasize the need for an efficient algorithm that quickly gets close to a solution: While for the H/ABM the inexact gradient descent algorithm is a viable option, it is not for GERDA. Here, the multilevel optimization approach is more efficient, although it also requires approximately 1.2~million GERDA simulations after the third iteration step. However, it comes with a significant reduction from $J^\text{ABM}(u_0) \approx 1\,189$ to $J^\text{ABM}(u_3) \approx 157$ after three iterations.

In the later phase, the agreement of the Hessian matrices of the fine and coarse models dominates the asymptotic convergence. Table~\ref{tab:experimental-convergence-rates} summarizes theoretical and experimental convergence rates for different coarse models acting as preconditioner. For the case $n_\mathrm{u} = 2$, the theoretical convergence rates were obtained by computing the condition numbers $\kappa$ of the Hessian matrices $H$ of a two-dimensional polynomial fit of degree five of the objective $J^\text{ABM}$ over the entire feasible optimization domain $U$ for the respective models. For the case $n_\mathrm{u} = 4$, we used a local four-dimensional quadratic polynomial fit of $J^\text{ABM}$ in an $\varepsilon$-environment around $u^\star$. In all cases the preconditioned problem leads to a faster convergence. The results are confirmed by the convergence plots (cf.~also Figure~\ref{fig:convergence-mock} for the H/ABM).

\begin{table}[t]
\centering
\caption{Theoretical and experimental convergence rates for the Multilevel Optimization Algorithm~\ref{alg:RMTR} for constant policies, i.e., $n_\mathrm{u} = 2$ and $T=1\,176$, and 25 days constant policies, i.e., $n_\mathrm{u} = 4$ and $T=1\,200$. The empirical convergence rate $\rho_\text{ODE, fine}$ is obtained using the ODE model as coarse-level and the H/ABM and GERDA, respectively, as fine-level models.}
\label{tab:experimental-convergence-rates}
\begin{tabular}{llccc}
\toprule
                              &                                            & \multicolumn{3}{c}{Fine model}                     \\ \cmidrule{3-5} 
                              &                                            & H/ABM    & GERDA    & H/ABM                        \\ \midrule
Dimension                     & $n_\mathrm{u}$                             & 2        & 2        & 4                            \\ \midrule
Condition                     & $\kappa(H_\text{ODE}^{-1}H_\text{fine})$   & 3.0300   & 116.8909 & 21.9801                      \\
                              & $\kappa(H_\text{H/ABM}^{-1}H_\text{fine})$ & --       & 7.0511   & --                           \\
                              & $\kappa(H_\text{fine})$                    & 35.1841  & 235.8803 & 131.8860                     \\ \midrule
Convergence rate              & $\rho_\text{ODE}$                          & 0.5037   & 0.9830   & 0.9130                       \\
                              & $\rho_\text{H/ABM}$                        & --       & 0.7516   & --                           \\
                              & $\rho_\text{fine}$                         & 0.9447   & 0.9916   & 0.9850                       \\ \midrule
Experimental convergence rate & $\rho_\text{ODE, fine}$                    & 0.7433   & 0.1934   & 0.6854                       \\ \bottomrule
\end{tabular}
\end{table}

%==============================================================================================

\section{Conclusion}
\label{sec:conclusion}
 
We presented a heterogeneous multilevel optimization approach combining a fine-level ABM with a coarse-level ODE to find exemplary non-phar\-ma\-ceu\-ti\-cal interventions in epidemic policy design. We compared this method with state-of-the-art algorithms applicable to our optimization problem. The multilevel optimization algorithm is expected to be faster with respect to the number of iterations as it (i) captures the nonlinear structures better than a first- order second-order Taylor model in the early iteration, see Section~\ref{sec:ml-opt}, and (ii) a well-matched coarse model serves as preconditioner and thus reduces ``zig-zagging'' slowing down gradient methods, see Table~\ref{tab:experimental-convergence-rates}. 

However, the theoretical speedup of the multilevel optimization algorithm compared to the inexact gradient descent algorithm is mostly compensated by the exploding number of samples needed for a high confidence in the descent direction the closer the current iterate gets to the true minimizer, see Figure~\ref{fig:barchart}. This is also the main bottleneck for the inexact gradient descent algorithm. In the initial phase, though, the superiority of the multilevel optimization algorithm is apparent, as it leads to solutions that may already be sufficient in a real-world scenario, since the mathematical optimum need not or cannot be reached. 

Additionally, we showed that the optimal policies obtained using directly the ABM differ drastically from one using the the stylized SIR-type ODE model or H/ABM (see Figures~\ref{fig:energylandscapes} or~\ref{fig:GERDA_opt_ctrl_cmp}), which suggests that using detailed ABMs directly for the design of optimal policies is necessary for good results. The above optimization framework can be applied not only to epidemiological models, for finding optimal phar\-ma\-ceu\-ti\-cal and non-phar\-ma\-ceu\-ti\-cal interventions, but to any ABM where appropriate and meaningful controls can be found and applied and suitable reduced models are available. Future research will address the efficient gradient approximation of stochastic dynamical systems such as ABMs.

\subsection*{Data Availability}
The data that support the findings of this study are generated based on the models given in Section~\ref{sec:epidemic-models} and do not use any external data. The algorithms are implemented in MATLAB and Python and can be found at \url{https://github.com/Henningston/MLoptABM}.

%==============================================================================================

\bibliographystyle{plain}
\bibliography{references}

\pagebreak

\clearpage
\appendix
\section*{Supplementary Material}
\subsection*{Stochastic Gradient Descent}

Gradient averaging can be done implicitly on the fly by taking sufficiently small steps in the direction of a single or average of a few samples. This is known as finite difference based stochastic approximation or Kiefer--Wolfowitz algorithm~\cite{KieferWolfowitz1952}. For the purpose of illustration, we give a rather basic scheme shown in Algorithm~\ref{alg:SGD}, where $\Pi_U$ denotes a projection to the admissible set $U$. With a decreasing step size $c/k$ for some fixed $c>0$, the iteration converges almost surely to a local minimizer $u_*$. Note that a sample size $n>1$ for gradient estimation, known as \emph{mini-batch}, leads to faster convergence in terms of iterations, though not sample number, and allows for embarrassingly parallel evaluation of $\nabla\EE[J^{\rm ABM}(u_k)]$, which makes larger sample sizes attractive on current computer systems. With an optimal balance of sampling and discretization error, i.e., $h=\mathcal{O}(k^{-1/6})$, a convergence rate of $\mathcal{O}(k^{-1/3})$ is obtained. Using improvements such as iterate averaging as in Algorithm~\ref{alg:SGD}, an improved convergence rate of $\mathcal{O}(k^{-1/2})$ can be obtained. Nevertheless, even with these improvements the convergence is sublinear, such that many samples are required.

\begin{algorithm}
    \textbf{Input} $c>0$, $u_0\in U$, $n\ge 1$, $K\ge 1$, $p<0$, $h_0 > 0$ \\
    \For{$k = 0, \dots, K$}{
        $h_k = k^p h_0$\\
        $u_{k} = \Pi_U\left(u_{k-1} - \frac{c}{k} \nabla_{h_k} \EE_n[J^{\rm ABM}(u_{k-1})]\right)$
        } 
    \textbf{Output} $u_* = \frac{1}{K}\sum_{k=1}^K u_k$
    \caption{Basic mini-batch Kiefer--Wolfowitz algorithm with iterate averaging}
    \label{alg:SGD}
\end{algorithm}

\end{document}